\documentclass{amsart}

\usepackage[english]{babel}
\usepackage[T1]{fontenc}
\usepackage{fullpage}
\usepackage[dvips]{epsfig}
\usepackage{amsfonts, amssymb, amsmath, stmaryrd, amsthm, euscript, mathrsfs}
\usepackage{color}
\usepackage{graphicx}
\usepackage{eurosym}
\usepackage{latexsym}
\usepackage{ifthen}
\usepackage{multirow}
\usepackage{url}
\usepackage{subfig}
\usepackage{booktabs}

\RequirePackage[colorlinks=true,citecolor=blue,urlcolor=blue]{hyperref}

\newtheorem{thm}{Theorem}

\newtheorem{prop}{Proposition}
\newtheorem{assumption}{Assumption}
\newtheorem{lem}{Lemma}

\newtheorem{rk}{Remark}

\newcommand{\poubelle}[1]{}

\newcommand{\E}{\mathbb{E}}

\renewcommand{\P}{\mathbb{P}}
\newcommand{\R}{\mathbb{R}}
\newcommand{\Z}{\mathbb{Z}}

\newcommand{\cF}{{\mathcal F}}

\newcommand{\cI}{{\mathcal I}}
\newcommand{\cJ}{{\mathcal J}}

\newcommand{\cM}{{\mathcal M}}
\newcommand{\cN}{{\mathcal N}}
\newcommand{\cP}{{\mathcal P}}

\newcommand{\cS}{{\mathcal S}}
\newcommand{\cT}{{\mathcal T}}

\newcommand{\indiq}{{{\bf 1}}}

\newcommand{\sigmaS}[1]{\ifthenelse{\equal{#1}{}}{\sigma}{\sigma_{#1}}}
\newcommand{\sigmaL}[1]{\ifthenelse{\equal{#1}{}}{\overline{\sigma}}{\overline{\sigma}_{#1}}}
\newcommand{\drift}[2]{b_{#1}^{#2}}
\newcommand{\Ws}[1]{\ifthenelse{\equal{#1}{}}{B}{B_{#1}}}
\newcommand{\Wl}[1]{\ifthenelse{\equal{#1}{}}{\overline{B}}{\overline{B}_{#1}}}

\newcommand{\etaS}[1]{\ifthenelse{\equal{#1}{}}{\eta}{\eta_{#1}}}
\newcommand{\etaL}[1]{\ifthenelse{\equal{#1}{}}{\overline{\eta}}{\overline{\eta}_{#1}}}

\newcommand{\diffE}[3]{\mathfrak{e}_{#1,#2}(#3)}

\definecolor{vertRemarque}{rgb}{.2, .9, .1}

\newcommand{\beqs}{\begin{eqnarray*}}
\newcommand{\enqs}{\end{eqnarray*}}
\newcommand{\beq}{\begin{eqnarray}}
\newcommand{\enq}{\end{eqnarray}}

\begin{document}

\title{Efficient volatility estimation in a two-factor model}

\author{Olivier~\textsc{F\'eron} \and Pierre \textsc{Gruet} \and Marc \textsc{Hoffmann}}

\address{Olivier F\'eron, EDF Lab Paris-Saclay and FiME, Laboratoire de Finance des March\'es de l'Energie, 91120 Palaiseau, France.}

\email{olivier-2.feron@edf.fr}

\address{Pierre Gruet, EDF Lab Paris-Saclay and FiME, Laboratoire de Finance des March\'es de l'Energie, 91120 Palaiseau, France.}

\email{pierre.gruet@edf.fr }

\address{Marc Hoffmann, Universit\'e Paris-Dauphine, CEREMADE, 75016 Paris, France}

\email{hoffmann@ceremade.dauphine.fr}

\begin{abstract}
We statistically analyse a multivariate HJM diffusion model with stochastic volatility. 
The  volatility process of the first factor is left totally unspecified while the volatility of the second factor is the product of an unknown process and an exponential function of time to maturity. This exponential term includes some real parameter measuring the rate of increase of the second factor as time goes to maturity. From historical data, we efficiently estimate the time to maturity parameter in the sense of constructing an estimator that achieves an optimal information bound in a semiparametric setting. We also identify nonparametrically the paths of the volatility processes and achieve minimax bounds. We address the problem of degeneracy that occurs when the dimension of the process is greater than two, and give in particular optimal limit theorems under suitable regularity assumptions on the drift process. We consistently analyse the numerical behaviour of our estimators on simulated and real datasets of prices of forward contracts on electricity markets.
\end{abstract}

\maketitle



\section{Introduction}

\subsection{Motivation and setting} \label{subsection:motivationSetting}

We address statistical estimation for multidimensional diffusion processes from historical data, with a volatility structure including both a parametric and a nonparametric components. We aim at achieving efficient estimation of a scalar parameter in the volatility, in presence of nonparametric nuisance, while providing point estimates of nonparametric components simultaneously. The processes of interest follow the multiple Brownian factor representation, as in the Heath-Jarrow-Morton (HJM) framework for forward rates, for instance in Heath \emph{et al.}~\cite{bib:HJM}, or for electricity forward contracts in Benth and Koekebakker~\cite{bib:BenthKoekebakkerStochasticModeling}.\\

Our setting is motivated by the context of prices of specific forward contracts, which are available on the electricity market. Interest rate models have been applied to the pricing of such contracts: see for instance Hinz \emph{et al.}~\cite{bib:hinz2005}, in which an analogy between interest rate models and forward contracts prices models is performed, the maturity in the former framework being a date of delivery in the latter. The factorial representation of the HJM framework has been precisely studied in Benth and Koekebakker~\cite{bib:BenthKoekebakkerStochasticModeling} to model the electricity forward curve, giving constraints in the volatility terms to ensure no arbitrage. Koekebakker and Ollmar~\cite{bib:forwardCurveDynamicsNordicElectricityMarket} perform a Principal Component Analysis to point out that two factors can explain 75\% of the electricity forward contracts in the Norwegian market, and more than 10 factors are needed to explain 95\%. They argue that, due to the non-storability of electricity, there is a weak correlation between short-term and long-term events. In Keppo \emph{et al.}~\cite{bib:KAHV}, a one-factor model is designed for each maturity date, having correlations between the Brownian motions for distinct dates. In Kiesel \emph{et al.}~\cite{bib:kiesel2009}, a two-factor model is described, with a specification of the volatility terms allowing to reproduce the classical behaviour of prices, especially the empirical evidence of the Samuelson effect (the volatility of prices increases as time to maturity decreases) and to ensure non-zero volatility for long-term forward prices.\par \medskip

On some filtered probability space $(\Omega, \cF, (\cF_t)_{t \geq 0}, \P)$, we consider a $d$-dimensional It\^o semimartingale $X = (X_t)_{t \geq 0}$ with components $X^j$, for $j=1,\ldots, d$, of the form
\begin{equation} 
        X_t^j = X_0^j + \int_0^t \drift{s}{j} ds + \int_0^t e^{-\vartheta (T_j-s)} \sigmaS{s} d\Ws{s} + \int_0^t \sigmaL{s} d\Wl{s}, \label{eqn:dynamiqueXti}
\end{equation}
where $X_0^j \in \R$ is an initial condition, $\Ws{} = (\Ws{t})_{t\geq 0}$ and $\Wl{}=(\Wl{t})_{t \geq 0}$ are two independent Brownian motions, $\vartheta$ and $T_j$ are positive numbers and $\sigmaS{} = (\sigmaS{t})_{t \geq 0}$, $\sigmaL{} = (\sigmaL{t})_{t \geq 0}$, $\drift{}{j} = (\drift{t}{j})_{t \geq 0}$ are c\`adl\`ag adapted processes. To avoid trivial situations, we assume that for some $T>0$, we have
\[
        T \leq T_1 < \ldots < T_d
\]
and that the $T_j$ are known. Moreover, we observe $X$ at times
\[
        0, \Delta_n, 2\Delta_n, \ldots, n\Delta_n=T.
\]
Asymptotics are taken as $n \rightarrow \infty$ with fixed $T$ once for all and throughout the paper. In this setting, it is impossible to identify the components $\drift{}{i}$, so  we are left with trying to estimate the parameter $\vartheta$ and the random components $t \leadsto \sigmaS{t}$ (or rather $\sigmaS{t}^2$) and $t \leadsto \sigmaL{t}$ (or $\sigmaL{t}^2$) over the time interval $[0,T]$ with the best possible rate of convergence. This is not always possible and will require regularity assumptions.\\

The statistical estimation of the volatility of a diffusion process observed over some period $[0,T]$ has long been studied for asymptotic regimes in which observation times asymptotically recover the whole observation period.  This carries over to the setting, considered here, where the unknown volatility -- as a parameter -- is random w.r.t. the filtration generated by the observation itself, see for instance \cite{bib:genonCatalotJacod1993, bib:GenonCatalotLaredoPicardOndelettes, bib:HoffmannIrregularSamplings, H4, HMS} and the references therein for a comprehensive study in both parametric and nonparametric settings. Concerning estimating a functional of the trajectory of the diffusion process
the chapters of Mykland and Zhang~\cite{bib:econometricsOfHFDataInStatisticalMethodsForSDE} and Jacod~\cite{bib:statisticsAndHFDataInStatisticalMethodsForSDE} present the most advanced problems related to the estimation of diffusion processes, together with important estimation results, stated in a general way. As integrated volatility can be estimated with the usual $\Delta_n^{-1/2}$-rate of convergence, the quality of its estimators may be assessed by looking at the limit law that one can get when writing a central limit theorem, and by looking at a minimal variance in some sense (usually, the limiting distribution is a mixture of a centred Gaussian variable, with random variance). Cl\'ement \emph{et al.}~\cite{bib:ClementDelattreGloterEfficiency} estimated some functionals of the volatility; in the diffusion model that they introduced, they prove an extension of H\'ajek convolution theorem, and are able to define some notion of efficiency, which is somehow related to our setting. The present paper is in line with these results from a methodological point of view.  However, we face several new structural aspects and questions that require a novel treatment and that lead to new results in some cases.
\begin{enumerate}
\item The dimension $d$ dictates in some sense the underlying regularity of the model. The (componentwise) quadratic variation identifies the parameters of the model asymptotically via its integrated volatility by the (componentwise) convergence 
$$
        \sum_{i=1}^{\lfloor \Delta_n^{-1} t \rfloor} (\Delta_i^n X^j)^2 \to \int_0^t
        \big(e^{-2\vartheta(T_j-t)}\sigmaS{t}^2+\sigmaL{t}^2\big)dt,\;\;j=1, ,\ldots, d
$$
        in probability as $\Delta_n\rightarrow 0$, valid for every $t \in [0,T]$ with the notation $\Delta_i^n X= X_{i\Delta_n}-X_{(i-1)\Delta_n}$ (componentwise). We thus identify the (random) function
        $$t \leadsto  e^{-2\vartheta(T_j-t)}\sigmaS{t}^2+\sigmaL{t}^2,\;\;j=1, ,\ldots, d$$
for $t \in [0,T]$.
       Since $\vartheta$, $\sigmaS{}$ and $\sigmaL{}$ are unknown, we readily see that the model cannot be identified for $d=1$. For $d=2$, we asymptotically recover the function
        $$t \leadsto 
        \left\{
        \begin{array}{l}
         e^{-2\vartheta(T_1-t)}\sigmaS{t}^2+\sigmaL{t}^2 \\
         e^{-2\vartheta(T_2-t)}\sigmaS{t}^2+\sigmaL{t}^2
        \end{array}
        \right.$$
and it is not obvious that the three parameters can be identified. This is possible by considering specific linear and quadratic combinations of the components of  the vectors $\Delta_i^n X$, $i=1,\ldots, n$ and this yields to estimators of $\vartheta$ and $\sigmaS{}, \sigmaL{}$ that achieve usual rates of convergence. In particular, the case $d=3$ (and beyond) becomes somehow degenerate, since we (discretely) observe a 3-dimensional diffusion model driven by a two dimensional Brownian motion $(B,\overline{B})$. For applications, this had been reported by Jeffrey \emph{et al.}~\cite{bib:jlnp2001} in a similar context. In that case, the same kind of techniques lead to estimator of $\vartheta$ that achieve the fast rate $\Delta_n^{-1}$. 
\item It is interesting to note that it is possible in dimension $d=2$ to separate $\vartheta$ from the nuisance nonparametric components $\sigmaS{}$ and $\sigmaL{}$. However, the substitution estimators are only rate-optimal. We need to modify our estimator of $\vartheta$ and introduce a completely different technique in order to obtain an optimal asymptotic variance. This is done by using semiparametric efficiency theory see {\it e.g.} the classical textbook of van der Vaart \cite{bib:vanDerVaartAsymptoticStatistics} in a context where the nonparametric nuisance parameters is a bivariate random volatility process $t \leadsto (\sigmaS{}_t^2,\sigmaL{}_t^2)$, which is new to the best of our knowledge.
\item In particular, we need to exhibit nonparametric estimators of $\sigmaS{}$ and $\sigmaL{}$ to be used in the correction efficient estimator of $\vartheta$. While nonparametric estimation of the time-dependent diffusion coefficient is classical, we do not have a completely usual nonparametric problem, since
\begin{enumerate}
        \item $\sigmaS{}$ and $\sigmaL{}$ are random themselves, so that we do not estimate them pointwise, instead we estimate the trajectories $t \leadsto (\sigmaS{}^2(\omega)_t,\sigmaL{}^2(\omega)_t)$ pointwise, as realisations of the volatility process;
        \item an increment $\Delta_i^n X$ is the sum of two stochastic integrals, in which the volatility processes have different regularities.
\end{enumerate}              
\end{enumerate}

\subsection{Main results and organisation of the paper}

In Section~\ref{subsection:estimationTheta}, we provide an estimator of $\vartheta$, based on quadratic variation, in the above observation scheme. While we cannot perform estimation when the number of observed processes $d$ is equal to $1$, the case $d=2$ is statistically regular, and by approaching the quadratic variation of $X^1$, $X^2$ and $X^2-X^1$, we derive an estimator $\widehat{\vartheta}_{2,n}$ of $\vartheta$, which is $\Delta_n^{-1/2}$-consistent. Using the theory of statistics for diffusion processes and relying on the tools of stable convergence in law, which are for instance summarized in \cite{bib:statisticsAndHFDataInStatisticalMethodsForSDE,bib:econometricsOfHFDataInStatisticalMethodsForSDE}, we show that
\[
        \Delta_n^{-1/2}(\widehat{\vartheta}_{2,n}-\vartheta)\to \cN(0,V_\vartheta(\sigmaS{},\sigmaL{})),
\]
stably in law, where, conditional on $\mathcal F_T$, the random variable  $\cN(0,V_\vartheta(\sigmaS{},\sigmaL{}))$ is centred Gaussian, with conditional variance $V_\vartheta(\sigmaS{},\sigmaL{})$, and possibly
 defined on an extension of the original probability space. When $d = 3$ (and beyond) the model is somehow degenerate because the $3$ marginal components of the process are driven by $2$ Brownian motions. The remaining source of randomness is the drift process, and we are able to construct a $\Delta_n^{-1}$-consistent estimator $\widehat{\vartheta}_{3,n}$ for $\vartheta$. This is possible as soon as $\drift{}{}$ has some integrated regularity in expectation, reminiscent of the so-called Besov regularity, as will be made precise by Assumption~\ref{asm:regulariteFonctions}. This enables us to obtain a satisfying limit theorem for $\Delta_n^{-1}(\widehat{\vartheta}_{3,n}-\vartheta)$, namely the convergence in probability to some $\cF$-measurable random variable. All our results in dimension $d=2,3$ are stated in Theorem~\ref{thm:loiAsymptotiqueHatAMoinsA}. The estimation of $\vartheta$ in $d=2$ however misses the optimal variance and we next turn to modifying $\widehat{\vartheta}_{3,n}$ in order to achieve the best possible variance. 
 
 In Section~\ref{subsection:estimationSigma}, we turn to the problem of estimating $\vartheta$ optimally when $d=2$. As a preliminary technical result, We first perform a relatively classical nonparametric estimation procedure to get point estimates of $\sigmaS{t}^2$ and $\sigmaL{t}^2$ when $d=2$ that will serve our later purposes of estimating $\vartheta$ in an optimal way.
We have to separate, in some way, the parts of the random increments that are linked to each of the Brownian integrals, to be able to get estimates of each process. We derive estimators $(\widehat{\sigmaS{}}_n^2,\widehat{\sigmaL{}}_n^2)$ of $(\sigmaS{}^2,\sigmaL{}^2)$. If the volatility processes are H\"older in expectation (Assumption~\ref{asm:sigmasHolder}) we show in Proposition~\ref{thm:estimationNPunePeriode} the convergence of the estimators with rate $\Delta_n^{-\alpha/(2\alpha+1)}$-consistent, where $\alpha$ is the lowest of two values of the H\"older regularities of $\sigmaS{}^2$ and $\sigmaL{}^2$. In Section~\ref{subsection:estimationEfficaceTheta}, relying to the theory of asymptotically efficient semiparametric estimation (for instance \cite{bib:vanDerVaartAsymptoticStatistics}) we compute a lower bound $V^{\mathrm{opt}}_\vartheta(\sigmaS{}, \sigmaL{})$ for the limit variance while estimating $\vartheta$ with $d=2$ observed processes, for deterministic volatility functions, in Theorem~\ref{thm:borneInferieureVarianceEstimationTheta}. As soon as $\sigmaL{}$ is not constant, this bound is lower than $V_\vartheta(\sigmaS{},\sigmaL{})$. We subsequently derive an estimator $\widetilde{\vartheta}_{2,n}$ such that
\[
        \Delta_n^{-1/2}(\widetilde{\vartheta}_{2,n}-\vartheta)\to\cN(0,V^{\mathrm{opt}}_\vartheta(\sigmaS{}, \sigmaL{}))
\]
stably in law, where conditional on $\mathcal F$, the random variable  $\cN(0,V^{\mathrm{opt}}_\vartheta(\sigmaS{}, \sigmaL{}))$ is centred Gaussian  with conditional variance $V^{\mathrm{opt}}_\vartheta(\sigmaS{},\sigmaL{})$. This is the main result of the paper. The estimator $\widetilde{\vartheta}_{2,n}$ is built upon the preliminary estimators $\widehat{\vartheta}_{2,n}$ and $(\widehat{\sigmaS{}}_n^2,\widehat{\sigmaL{}}_n^2)$ of $(\sigmaS{}^2,\sigmaL{}^2)$. Moreover, it is asymptotically efficient in the sense that it achieves the minimal conditional variance $V^{\mathrm{opt}}_\vartheta(\sigmaS{},\sigmaL{})$ among all possible $\Delta_n^{-1/2}$-consistent estimators that are asymptotically centred mixed normal.
We perform numerical experiments in Section~\ref{section:resultatsNumeriques}, using both simulated and real data from the electricity forward markets in order to compare the behaviours of the estimators in various configurations. The proofs are delayed until Section~\ref{section:preuves}.

\section{Construction of the estimators and convergence results}
\subsection{Rate-optimal estimation of $\vartheta$} \label{subsection:estimationTheta}

In this section, we build a (preliminary) estimator of $\vartheta$, denoted by $\widehat{\vartheta}_{d,n}$ in each of the situation $d=2$ or $d=3$.
\subsubsection*{The case $d=1$}
In that setting, it is impossible to identify $\vartheta$ from data $X_{i\Delta_n}, i=1,\ldots, n$ asymptotically when $t\leadsto \sigmaS{t}$ and $t \leadsto \sigmaL{t}$ are unknown. Indeed $X$ has the same law under the choice of $(\vartheta,\sigmaS{},\sigmaL{})$ and $(\vartheta+1, e^{T_1-\cdot}\sigmaS{}, \sigmaL{})$.
\subsubsection*{The case $d=2$}
This is the statistically most regular case. Set, as usual $\Delta_i^n X= X_{i\Delta_n}-X_{(i-1)\Delta_n}$ (componentwise).
From the convergences
\[
        \sum_{i=1}^n (\Delta_i^n X^j)^2 \to \int_0^T\big(e^{-2\vartheta(T_j-t)}\sigmaS{t}^2+\sigmaL{t}^2\big)dt,\;\;j=1,2
\]
and
\[
        \sum_{i=1}^n (\Delta_i^n X^2-\Delta_i^n X^1)^2 \to \int_0^T(e^{-\vartheta T_2}-e^{-\vartheta T_1})^2e^{2\vartheta t}\sigmaS{t}^2dt
\]
in probability, we also obtain the convergence of the ratio
\[
        \Psi_{T_1,T_2}^n =  \frac{\sum_{i=1}^n (\Delta_i^n X^2-\Delta_i^n X^1)^2}{\sum_{i=1}^n \big((\Delta_i^n X^2)^2-(\Delta_i^n X^1)^2\big)} \rightarrow 
\frac{\big(e^{-\vartheta T_2}-e^{-\vartheta T_1}\big)^2}{e^{-2\vartheta T_2}-e^{-2\vartheta T_1}}
= \psi_{T_1,T_2}(\vartheta),
\]
in probability. The function $\vartheta \leadsto \psi_{T_1,T_2}(\vartheta)$ maps $(0,\infty)$ onto $(-1,0)$ and this leads to a first 
estimation strategy by setting 
\[
        \widehat{\vartheta}_{2,n} = \psi_{T_1,T_2}^{-1}\big(\Psi_{T_1,T_2}^{n}\big)
\]
whenever $\Psi_{T_1,T_2}^{n} \in (-1,0)$ and $0$ otherwise.

\subsubsection*{The case $d=3$}

Since $X$ is driven by two Brownian motions, the underlying statistical model becomes degenerate. 
Indeed, assume first that $\drift{}{1}=\drift{}{2}=\drift{}{3}$. Then, we readily obtain 
\[
        \frac{\Delta_i^n X^2 - \Delta_i^n X^1}{\Delta_i^n X^3 - \Delta_i^n X^2}=\frac{e^{-\vartheta T_2}-e^{-\vartheta T_1}}{e^{-\vartheta T_3}-e^{-\vartheta T_2}}
\]
which is invertible as a function of $\vartheta$.
It is thus possible to identify $\vartheta$ exactly from the observation of a single increment of $X$ itself! When the $\drift{}{j}$ are not all equal, the situation is still somehow degenerate, as we can eliminate all volatility components by taking linear combinations of the observed increments. The lowest-order remaining term is the drift process, with increments that are of smaller order than that of the diffusion part. Therefore, we may hope to exhibit estimators of $\vartheta$ that converge at a faster rate then $\Delta_n^{1/2}$. More specifically, we have the convergence
\[
        \Psi_{T_1,T_2, T_3}^{n}  = \frac{\sum_{i=1}^n(\Delta_i^n X^3-\Delta_i^n X^2)^2}{\sum_{i=1}^n(\Delta_i^n X^2-\Delta_i^n X^1)^2} \to \Big(\frac{e^{-\vartheta T_3}-e^{-\vartheta T_2}}{e^{-\vartheta T_2}-e^{-\vartheta T_1}}\Big)^2 = \psi_{T_1,T_2,T_3}(\vartheta),
\]
say, in probability. The function $\vartheta \leadsto  \psi_{T_1,T_2,T_3}(\vartheta)$ maps $(0,\infty)$ onto $\big(0, \big(\frac{T_3-T_2}{T_2-T_1}\big)^2\big)$ and it can be checked with elementary calculus arguments that it is also invertible,
leading to the estimator
\[
        \widehat \vartheta_{3,n} = \psi_{T_1,T_2,T_3}^{-1}\big(\Psi_{T_1,T_2, T_3}^{n}\big)
\]
whenever $\Psi_{T_1,T_2, T_3}^n \in \big(0, \big(\frac{T_3-T_2}{T_2-T_1}\big)^2\big)$ and $0$ otherwise.

\subsubsection*{Convergence results}

Remember that $T$ is fixed and thus asymptotics are taken as $n \rightarrow \infty$ or equivalently $\Delta_n \rightarrow 0$. We need some assumption about the regularity of the processes $b$, $\sigmaS{}$ and $\sigmaL{}$. For a random process $X=(X_t)_{0 \leq t \leq T}$, introduce the following modulus of continuity:
\[
        \omega(X)_t=\sup_{|h| \leq t}\Big(\int_0^T \E\big[(X_{s+h}-X_s)^2\big]ds\Big)^{1/2}.
\]

\begin{assumption} \label{asm:regulariteFonctions} For some constant $c_{\min} >0$, we have $\inf_{t \in [0,T]}\min(\sigmaS{}_t,\sigmaL{}_t) \geq c_{\min}$. 
Moreover, for some $s>1/2$, we have $\sup_{t \in [0,T]} t^{-s}\omega(\drift{}{j})_t <\infty$ for every $j=1,\ldots, d$. 
\end{assumption}
 
\begin{rk} \label{rk: besov restrict}
Assumption \ref{asm:regulariteFonctions} is technical. The first condition on the ellipticity of the volatility processes enables one to apply Girsanov theorem to get rid of the drift terms $b^j$ in some cases. The second assumption is related to the pathwise smoothness property of the drift $b^j$ in $L^2$ and is reminiscent of a Besov regularity. For instance, a sufficient condition is that $t \leadsto b^j_t$ is almost-surely $s$-H\"older continuous. In particular, we need a slightly stronger smoothness than that of typical Brownian paths.     
\end{rk}

To state the convergence results, we need some notation. Set
\[
        \overline{b}_t=2(e^{-\vartheta T_2}-e^{-\vartheta T_1})(e^{-\vartheta T_3}-e^{-\vartheta T_2})\big((e^{-\vartheta T_2}-e^{-\vartheta T_1})(\drift{t}{3}-\drift{t}{2})-(e^{-\vartheta T_3}-e^{-\vartheta T_2})(\drift{t}{2}-\drift{t}{1})\big)
\]
and
\[
        \widetilde{b}_T=(e^{-\vartheta T_2}-e^{-\vartheta T_1})^2\int_0^T(\drift{t}{3}-\drift{t}{2})^2dt-(e^{-\vartheta T_3}-e^{-\vartheta T_2})^2\int_0^T(\drift{t}{2}-\drift{t}{1})^2dt.
\]
Define also
\[
        D_3=(e^{-\vartheta T_3}-e^{-\vartheta T_2})\big((e^{-\vartheta T_3}-e^{-\vartheta T_2})(T_2 e^{-\vartheta T_2}-T_1e^{-\vartheta T_1})-(e^{-\vartheta T_2}-e^{-\vartheta T_1})(T_3e^{-\vartheta T_3}-T_2e^{-\vartheta T_2})\big).
\]
\begin{thm} \label{thm:loiAsymptotiqueHatAMoinsA}
        Work under Assumption~\ref{asm:regulariteFonctions}.
        \begin{enumerate}
                \item For $d=2$, we have 	      
                \[
		       \Delta_n^{-1/2}(\widehat{\vartheta}_{2,n}-\vartheta) \rightarrow  \cN\big(0,V_\vartheta(\sigmaS{}, \sigmaL{})\big)
                \]
		in distribution as $n\rightarrow \infty$, where $\mathcal N\big(0,V_\vartheta(\sigmaS{}, \sigmaL{})\big)$ is a random variable which, conditionally on $\cF_T$, is centred Gaussian with variance 
		\[
                        V_\vartheta(\sigmaS{}, \sigmaL{}\big)=\frac{1}{(T_2-T_1)^2}(e^{\vartheta T_2}-e^{\vartheta T_1})^2 \frac{\int_0^T e^{2\vartheta t}\sigmaS{t}^2\sigmaL{t}^2 dt}{\big(\int_0^Te^{2\vartheta t}\sigmaS{t}^2 dt\big)^2}.
                \]
      \item Assume moreover that $b^j=0$. For $d=3$ we have
                \[
                     \Delta_n^{-1}(\widehat{\vartheta}_{3,n}-\vartheta) \rightarrow  \frac{\widetilde{b}_T+\int_0^T\bar{b}_te^{\vartheta t}\sigmaS{t}dB_t}{2(e^{-\vartheta T_2}-e^{-\vartheta T_1}) D_3\int_0^Te^{2\vartheta t}\sigmaS{t}^2dt} 
                \]
                in probability as $n\rightarrow \infty$.
        \end{enumerate}
\end{thm}

\begin{rk}
In the proof in Section \ref{proof of thm 1} below, we obtain $(1)$ when $d=2$ for $b^j=0$, $j=1,2$. The general case is obtained via a change of probability argument thanks to Girsanov theorem that can be applied under Assumption \ref{asm:regulariteFonctions}. This can no longer be applied when $d=3$ since the 3-dimensional semimartingale $X^j$ relies on the 2-dimensional Brownian motion $(B, \overline{B})$ solely. We therefore obtain $(2)$ under the restriction $b^j=0$. The general case can however be obtained by specifically evaluating additional superoptimal error terms involving $b^j$ at a significant cost of technical length.  
\end{rk}

\subsection{Rate-optimal estimation of the volatility processes} \label{subsection:estimationSigma}

\subsubsection*{Construction of an estimator}

In order to improve $\widehat{\vartheta}_{2,n}$  and building an asymptotically efficient estimator of $\vartheta$, we shall build pointwise rate-optimal nonparametric estimators of $\sigmaS{}^2$ and $\sigmaL{}^2$.  We start with the observation that for any sufficiently regular test function $g:[0,T]\rightarrow\R$, we have, for any $j=1,\ldots, d$,
\begin{equation} \label{identify sigma}
        \sum_{i = 1}^n g((i-1)\Delta_n)\big(\Delta_i^n X^j\big)^2 \rightarrow \int_0^T g(s)d\langle X^j \rangle_s = \int_0^Tg(s)\big(e^{-2\vartheta(T_j-s)}\sigmaS{s}^2+\sigmaL{s}^2\big)ds
\end{equation}
in probability as $n \rightarrow \infty$. Therefore, picking a function $g$ that mimics a Dirac mass at $t$, we can asymptotically identify 
\[
        e^{-2\vartheta(T_1-t)}\sigmaS{t}^2+\sigmaL{t}^2 \;\;\text{and}\;\; e^{-2\vartheta(T_2-t)}\sigmaS{t}^2+\sigmaL{t}^2
\]
by applying \eqref{identify sigma} for $j=1, 2$ for a sequence $g_n$ that converges to $\delta_t$ weakly. We thus identify $\sigmaS{t}^2$ and $\sigmaL{t}^2$ as well by inverting a $2\times 2$ linear system, namely
\[
	\begin{pmatrix} \sigmaS{t}^2 \\ \sigmaL{t}^2\end{pmatrix}=\cM(\vartheta)_t\begin{pmatrix}e^{-2\vartheta(T_1-t)}\sigmaS{t}^2+\sigmaL{t}^2 \\ e^{-2\vartheta(T_2-t)}\sigmaS{t}^2+\sigmaL{t}^2\end{pmatrix}
\]
where
\[
	\cM(\vartheta)_t=\frac{1}{e^{-2\vartheta(T_1-t)}-e^{-2\vartheta(T_2-t)}}\begin{pmatrix} 1 & -1 \\ -e^{-2\vartheta(T_2-t)} & e^{-2\vartheta(T_1-t)}\end{pmatrix} \text{.}
\]
For a threshold $\varpi_n>0$ and a bandwidth $h_n>0$, define the estimators
\begin{equation} \label{eqn:estimateurSigmaCetSigmaL}
	\begin{pmatrix}
	\widehat{\sigmaS{}}_{n,t}^2 \\ \widehat{\sigmaL{}}_{n,t}^2
	\end{pmatrix}
	=
	h_n^{-1}\cM(\max(\widehat{\vartheta}_{2,n},\varpi_n))_t \sum_{t-h_n \leq (i-1)/n < t} \begin{pmatrix}(\Delta_i^n X^1)^2 \\ (\Delta_i^n X^2)^2
	\end{pmatrix}.
\end{equation}
The bandwidth $h_n$ is set below to balance both bias and variance, while $\varpi_n >0$ guarantees the well-posedness of the estimator. Moreover, the exponential functions appearing in the denominator of $\widehat{\sigmaS{}}_{n,t}^2$ and not in its numerator may imply an ill-posedness of the estimation problem. In particular, when $\widehat{\vartheta}_{2,n}$ is high, one can expect very high values for $\widehat{\sigmaS{}}_{n,t}^2$. This is what happens in the numerical results of Section~\ref{sec_res_simulation_data}.

\subsubsection*{Convergence result}

We need an additional regularity assumption on the volatility processes $\sigmaS{}$ and $\sigmaL{}$.

\begin{assumption} \label{asm:sigmasHolder}
        There exists $\alpha \geq 1/2$ and a constant $c_\alpha >0$ such that for every $t,s \in [0,T]$, we have
        \begin{equation} \label{reg sigma}
                \E\big[(\sigmaS{t}^2-\sigmaS{s}^2)^2\big] +\E\big[(\sigmaL{t}^2-\sigmaL{s}^2)^2\big]\leq c_{\alpha}|t-s|^{2\alpha}.
        \end{equation}
\end{assumption}

\begin{rk} Assumption \ref{asm:sigmasHolder} is technical and related to the pathwise smoothness property of the volatility processes $\sigmaS{}$ and $\sigmaL{}$ in $L^\infty$. Like Assumption \ref{asm:regulariteFonctions}, it is reminiscent of a Besov regularity. In particular, in Theorem \ref{thm:estimationEfficace}, we need $\alpha >1/2$, a slightly stronger smoothness than that of typical Brownian paths for which  $\alpha= 1/2$ exactly.     
\end{rk}
\begin{prop} \label{thm:estimationNPunePeriode}
        Work under Assumptions \ref{asm:regulariteFonctions} and  \ref{asm:sigmasHolder}. Specify
                $h_n = \Delta_n^{1/(2\alpha+1)}$ and
       let $\varpi_n \rightarrow 0$. 
        Then the sequences
        \[
                \Delta_n^{-\alpha/(2\alpha+1)}\big(\widehat{\sigmaS{}}_{n,t}^2-\sigmaS{t}^2)\;\;\text{and}\;\; \Delta_n^{-\alpha/(2\alpha+1)}(\widehat{\sigmaL{}}_{n,t}^2-\sigmaL{t}^2)
        \]
        are tight, uniformly in $t$ over compact sets included in $(0,T]$.
\end{prop}

\subsection{Efficient estimation of $\vartheta$ when $d=2$} \label{subsection:estimationEfficaceTheta}

We may now look for the best attainable variance among rate-optimal estimators of $\vartheta$ that are asymptotically Gaussian.
However, we do not have a statistical model in the classical sense, with parameters $(\vartheta,\sigmaS{},\sigmaL{})$ since $t \leadsto\sigmaS{}_t$ and $t \leadsto \sigmaL{}_t$ are random processes themselves. In order to bypass this difficulty, we first restrict our attention to the case where $\sigmaS{}$ and $\sigmaL{}$ are deterministic functions, which enables us to identify our data within a semiparametric regular statistical model. Thanks to classical bounds on semiparametric estimation, we can explicitly compute the optimal (best achievable) variance $V^{\text{opt}}_\vartheta(\sigmaS{}, \sigmaL{})$. In a second step, allowing $\sigmaS{}$ and $\sigmaL{}$ to be random again, we build a one-step correction of our preliminary estimator $\widehat{\vartheta}_{2,n}$ which has the property of being asymptotically mixed Gaussian, with (conditional) variance equal to $V^{\text{opt}}_\vartheta(\sigmaS{}, \sigmaL{})$, \emph{i.e.} thus achieving the optimal variance along deterministic paths.

\subsubsection*{Lower bounds}

Consider the statistical experiment $\mathcal{E}^n$ generated by data $(\Delta_i^nX^1, \Delta_i^nX^2,i=1,\ldots, n)$ with
\begin{equation} \label{diff gaussienne}
X_t^j = X_0^j + \int_0^t e^{-\vartheta (T_j-s)} \sigmaS{s} d\Ws{s} + \int_0^t \sigmaL{s} d\Wl{s},\;\;j=1,2,
\end{equation}
with parameter $(\vartheta, \sigmaS{}, \sigmaL{}) \in \Theta \times \Sigma(c,\tilde{c})$, with $\Theta = (0,\infty)$ and $\Sigma(c,\tilde{c})$ being the space of positive (deterministic) functions $(\sigmaS{}, \sigmaL{})$ defined on $[0,T]$, satisfying \eqref{reg sigma} of Assumption \ref{asm:sigmasHolder} with constant $c$ and satisfying moreover $\tilde c \leq \inf_t \sigma_t \leq \sup_t \sigma_t\leq c$ for some $\tilde{c}>0$. Note that \eqref{diff gaussienne} is a Gaussian diffusion obtained from \eqref{eqn:dynamiqueXti} when restricting $(\sigmaS{}, \sigmaL{})$ to deterministic functions and setting $b^j=0$. 

\begin{thm} \label{thm:borneInferieureVarianceEstimationTheta}
        Let $\widehat{\vartheta}_n$ be an estimator of $\vartheta$ in the experiment $\mathcal{E}^n$ such that $\Delta_n^{-1/2}(\widehat{\vartheta}_n- \vartheta)$ converges to $ \cN\big(0,V_\vartheta(\sigmaS{}, \sigmaL{})\big)$ in distribution as $n \rightarrow \infty$. Then
        \[
                V_\vartheta(\sigmaS{}, \sigmaL{}) \geq  V^{\mathrm{opt}}_\vartheta(\sigmaS{}, \sigmaL{}) =  \frac{1}{(T_2-T_1)^2}(e^{\vartheta T_2}-e^{\vartheta T_1})^2\Big(\int_0^T \frac{e^{2\vartheta t}\sigmaS{t}^2}{\sigmaL{t}^2} dt\Big)^{-1}.
        \]
\end{thm}

Using Cauchy-Schwarz inequality, it is easy to prove that the expression of the limit variance is equal to the one we got in Theorem~\ref{thm:loiAsymptotiqueHatAMoinsA} for $\Delta_n^{-1/2}(\widehat{\vartheta}_{2,n}-\vartheta)$ if and only if $\sigmaL{}$ is constant over the interval $[0,T]$. In the general case of an arbitrary volatility process, we have a suboptimal result. This is the motivation of out next step, toward asymptotic efficiency.

\subsubsection*{Construction of an efficient procedure}

This is the most delicate part of the paper. By representation \eqref{diff gaussienne}, we see that the $(\Delta_i^nX^1,\Delta_i^n X^2)$ are independent for $i=1,\ldots, n$. Moreover, $(\Delta_i^nX^1,\Delta_i^n X^2)$ is a centred Gaussian, with explicit covariance 
structure
\begin{align*}
  \E\big[(\Delta_i^n X^1)^2\big]  & = \int_{(i-1)\Delta_n}^{i\Delta_n}e^{-2\vartheta(T_1-t)}\sigmaS{t}^2dt + \int_{(i-1)\Delta_n}^{i\Delta_n}\sigmaL{t}^2dt,\\
        \E\big[(\Delta_i^n X^2)^2\big] & = \int_{(i-1)\Delta_n}^{i\Delta_n}e^{-2\vartheta(T_2-t)}\sigmaS{t}^2dt+\int_{(i-1)\Delta_n}^{i\Delta_n}\sigmaL{t}^2dt,\\ 
        \E\big[\Delta_i^n X^1 \Delta_i^nX^2\big] & = \int_{(i-1)\Delta_n}^{i\Delta_n}e^{-\vartheta(T_1+T_2-2t)}\sigmaS{t}^2dt+\int_{(i-1)\Delta_n}^{i\Delta_n}\sigmaL{t}^2dt.
\end{align*}
Let us further denote by $f_{\vartheta, \sigmaS{}, \sigmaL{}}^{i,n}$ its density function w.r.t. the Lebesgue measure on $\R^2$.
If the nuisance parameters $(\sigmaS{},\sigmaL{})$ were known, then an optimal (efficient) procedure could be obtained by a one-step correction of the type
\[
        \widehat{\vartheta}_{n}=\widehat{\vartheta}_{2,n}+\frac{\sum_{i=1}^n \ell_{\vartheta = \widehat{\vartheta}_{2,n}, \sigmaS{}, \sigmaL{}}^i(\Delta_i^nX^1,\Delta_i^n X^2)}{\sum_{i=1}^n \big(\ell^i_{\vartheta = \widehat{\vartheta}_{2,n}, \sigmaS{}, \sigmaL{}}(\Delta_i^nX^1,\Delta_i^n X^2)\big)^2}
\]
where 
\begin{equation} \label{eq: def score}
\ell_{\vartheta, \sigmaS{}, \sigmaL{}}^{i,n}(\Delta_i^nX^1,\Delta_i^n X^2) = \partial_\vartheta \log f_{\vartheta, \sigmaS{}, \sigmaL{}}^{i,n}(\Delta_i^nX^1,\Delta_i^n X^2)
\end{equation} 
is the score function associated to $(\Delta_i^nX^1,\Delta_i^n X^2)$, see for instance Section~8.9 in \cite{bib:vanDerVaartAsymptoticStatistics}. However, this oracle procedure is not achievable and we need to invoke the theory of semiparametric efficiency (see for instance Chapter 25 of \cite{bib:vanDerVaartAsymptoticStatistics}). In the presence of an extra nuisance parameter $(\sigmaS{}, \sigmaL{})$, we consider instead the so-called efficient score
\[
        \widetilde{\ell}_{\vartheta, \sigmaS{}, \sigmaL{}}^{i,n} = \ell_{\vartheta, \sigmaS{}, \sigmaL{}}^{i,n}- \Pi \ell_{\vartheta, \sigmaS{}, \sigmaL{}}^{i,n},
\]
where $\Pi$ is the projection operator onto the tangent space associated to a one-dimensional perturbation around the true (unknown) value $(\sigmaS{}, \sigmaL{})$. It turns out that we indeed have a simple and explicit formula for $\widetilde{\ell}_{\vartheta, \sigmaS{}, \sigmaL{}}^{i,n}$ which enables us to derive a one-step correction formula using $\widetilde{\ell}_{\vartheta, \sigmaS{}, \sigmaL{}}^{i,n}$ and plug-in estimators in order to achieve the optimal bound.\\

For technical reason, we replace $\widehat{\vartheta}_{2,n}$ by $\Delta_n^{1/2}\lfloor \Delta_n^{-1/2}\widehat{\vartheta}_{2,n}\rfloor$ and we still write $\widehat{\vartheta}_{2,n}$ for simplicity. Likewise, we implicitly replace the estimators $\widehat{\sigmaL{}}_{n,t}^2$ defined in \eqref{eqn:estimateurSigmaCetSigmaL} by $\max\big(\min(\widehat{\sigmaL{}}_{n,t}^2, \tilde{c}^2),c\big)$, where $\tilde{c}$ is the lower bound associated to $\Sigma(c,\tilde{c})$ in the definition of the experiment $\mathcal{E}^n$. Set
\[
        \widetilde{\ell}_{\vartheta, \mu}^{i,n}(\Delta_i^nX^1,\Delta_i^n X^2)  = \frac{(\Delta_i^n X^2-\Delta_i^n X^1)(\Delta_i^n X^2-e^{-\vartheta (T_2-T_1)}\Delta_i^n X^1)e^{-\vartheta(T_2-T_1)}(T_2-T_1)}{(1-e^{-\vartheta (T_2-T_1)})^3 \Delta_n \mu}
\]
for $i=1,\ldots, n$ and $(\vartheta, \mu) \in (0,\infty) \times [\tilde c,c]$. We are ready to state the main result of the paper. We may assume here that the data $(\Delta_i^n X^j)$ are obtained under \eqref{eqn:dynamiqueXti} and therefore $b^j \neq 0$ here while the definition of $\widehat \vartheta_{n,3}$ remains unchanged.
\begin{thm} \label{thm:estimationEfficace}
        Work under Assumptions \ref{asm:regulariteFonctions} and \ref{asm:sigmasHolder} with $\alpha >1/2$. 
        \begin{enumerate}
        \item In the Gaussian experiment $\mathcal E^n$ generated by \eqref{diff gaussienne} the efficient score for the parameter $\vartheta$ associated to $(\Delta_i^nX^1,\Delta_i^n X^2)$ is given by   
        $$\widetilde{\ell}_{\vartheta, \Delta_n^{-1}\int_{(i-1)\Delta_n}^{i\Delta_n} \sigmaL{t}^2dt}^{i,n}(\Delta_i^nX^1,\Delta_i^n X^2).$$ 
\item   In the general experiment generated by \eqref{eqn:dynamiqueXti} allowing for $b^j \neq 0$, if moreover $\P\big((\sigmaS{}, \sigmaL{}) \in \Sigma(c,\tilde{c})\big)=1$ for some $c, \widetilde c>0$, the estimator $\widetilde{\vartheta}_{2,n}$ defined by
        \[
                \widetilde{\vartheta}_{2,n} = \widehat{\vartheta}_{2,n}+\frac{\sum_{i\in\cI_n} \widetilde{\ell}_{\widehat{\vartheta}_{2,n}, \widehat{\sigmaL{}}_{n,(i-1)\Delta_n}^2}^{i,n}(\Delta_i^nX^1,\Delta_i^n X^2)}{\sum_{i \in \cI_n} \big(\widetilde{\ell}_{\widehat{\vartheta}_{2,n}, \widehat{\sigmaL{}}_{n,(i-1)\Delta_n}^2}^{i,n}(\Delta_i^nX^1,\Delta_i^n X^2)\big)^2}
        \]
        with $\cI_n=\{i=1,\ldots,n, h_n \leq (i-1)\Delta_n < T\}$
        satisfies
        \[
                \Delta_n^{-1/2}\big(\widetilde{\vartheta}_{2,n}-\vartheta\big) \rightarrow \cN\big(0, V^{\mathrm{opt}}_\vartheta (\sigmaS{}, \sigmaL{})\big)
        \]
        in distribution as $n\rightarrow \infty$, where  $\cN\big(0, V^{\mathrm{opt}}_\vartheta (\sigmaS{}, \sigmaL{})\big)$ denotes a random variable that
        conditionally on $\cF_T$, centred Gaussian with (conditional) variance $V^{\mathrm{opt}}_\vartheta (\sigmaS{}, \sigmaL{})$.
        \end{enumerate}
\end{thm}
This result shows that the lower bound $V^{\mathrm{opt}}_\vartheta (\sigmaS{}, \sigmaL{})$ can be attained, and therefore that efficient estimation can be performed beyond deterministic volatility functions. This is reminiscent of the LAMN property when considering random volatilities in parametric estimation and can be interpreted as a semiparametric analog, although we do not have a complete theory at this stage. A recommended reference in that direction is Genon-Catalot and Jacod \cite{bib:genonCatalotJacod1993} where the same methodology is applied in a parametric context.

\section{Numerical implementation and real data analysis} \label{section:resultatsNumeriques}

\subsection{Electricity forward contracts}
The prices of existing forward contracts in the electricity markets are characterised by three time components: the quotation date $t$ and the dates $T_s$ and $T_e$ of respectively starting and ending power delivery. Therefore, a forward contract $F(t,T_s,T_e)$ will deliver to the holder 1~MWh of electricity continuously between dates $T_s$ and $T_e$. Such a contract may be bought during a quotation period $[t_0,T]$ with $T<T_s$ and it is no more available once $t>T$. Typical observed contracts are of various delivery periods: one week, one month, one quarter (three~months), one season (6~months) or one year. In this study we only consider the 6 observable monthly contracts (from 1 to 6 month-ahead) to estimate $\vartheta$ and the volatility processes $\sigmaS{}$ and $\sigmaL{}$. Also, for simplicity, we will drop $T_e$ from the notation. In the context of simulated data, we will simulate prices of $F(t,T_s)=F(t,T_s,T_e)$, the forward delivering continuously 1~MWh during the period $[T_s,T_e]$. In the context of real data, the price $F(t,T_s)$ is observable.

\subsection{Results on simulated data} \label{sec_res_simulation_data}

The objective of this section is to study the estimators' behaviour on a simulated data set, where the log-prices of the forward contracts are simulated according to the two-factor model described in \eqref{eqn:dynamiqueXti}. The parameter values are chosen to be close to values estimated on real data: in \cite{bib:kiesel2009}, the volatility processes are constant, and the estimated values are $\sigmaS{}=0.37$~yr$^{-1/2}$ and $\sigmaL{}=0.15$~yr$^{-1/2}$. Here we use a CIR-like model (the Cox-Ingersoll-Ross model for interest rates has been introduced in \cite{bib:CIRmodel}, in 1985), to emphasize the fact that our model may also be used in the context of interest rates modeling (this is indeed where it comes from, see \cite{bib:hinz2005}). Our parameters are
\[
        \drift{t}{j}=3.65\cdot 10^{-1}(\log(30)-X_t^j), \sigmaS{t}=0.37 \Sigma_t^d  \text{ and } \sigmaL{t}=0.15 \Sigma_t^d,
\]
with $\Sigma_t^d=\sqrt{\frac{1}{d}\sum_{j=1}^d X_t^j}$, which is the square root of the average of the $d$ quoted log-prices.
We adopt various values of $\vartheta$ (values in yr$^{-1}$): 1.4, 10, 40. The first value is the estimated parameter shown in \cite{bib:kiesel2009} and the others are chosen to cover a wide range of possible values to observe different behaviours of our estimators. Finally, the initial value of each simulated log-price series is the logarithm of a random variable taken uniformly over the interval $[20,40]$, which is an usual range for prices in the market of forward contracts on electricity (see also the constant $30$ in the drift, in the center of that interval).
We consider different simulation configurations, all related to the situations we are facing on real data. 
\begin{itemize}
	\item 2 processes (1 month-ahead and 2 month-ahead) observed on $n=100$ dates, with $T=T_1=150$ and $T_2=181$ days.
	\item 3 processes (1 month-ahead to 3 month-ahead) observed on $n=80$ dates, with $T=T_1=120$, $T_2=150$ and $T_3=181$ days.
	\item 4 processes (1 month-ahead to 4 month-ahead) observed on $n=60$ dates, with $T=T_1=90$, $T_2=120$, $T_3=151$ and $T_4=181$ days.
	\item 5 processes (1 month-ahead to 5 month-ahead) observed on $n=40$ dates, with $T=T_1=59$, $T_2=90$, $T_3=120$, $T_4=151$ and $T_5=181$ days.
	\item 6 processes (1 month-ahead to 6 month-ahead) observed on $n=20$ dates, with $T=T_1=31$, $T_2=59$, $T_3=90$, $T_4=120$, $T_5=151$ and $T_6=181$ days.
\end{itemize}
The decreasing number of observations corresponds to the configuration observed with real data: 2 monthly contracts with fixed delivery dates are jointly observed on working days during 5 months (around 100 quotation dates) whereas 6 monthly contracts can be jointly observed only during 1 month (around 20 quotation dates). The number of observations is a bit low, as we are relying on asymptotic results.\\
 
For each configuration, we perform 100,000 simulations. Recall that we denote by $\widehat{\vartheta}_{j,n}$ the estimator of $\vartheta$ from the configuration where $j$ processes are observed, and also by $\widetilde{\vartheta}_{2,n}$ the efficient estimator as described in Section~\ref{subsection:estimationEfficaceTheta}, available in the configuration of 2 observed processes. Although we have not proved that the estimator $\widetilde{\vartheta}_{2,n}$ is $\Delta_n^{-1/2}$-consistent and that it reaches the lower bound for the limit variance when $\alpha=1/2$, we have not got any numerical evidence against that possibility.
Tables~\ref{tab:simulation1}, \ref{tab:simulation2} 
and \ref{tab:simulation4} give the estimation results for $\vartheta=1.4$, $10$
and $40$~yr$^{-1}$, respectively. In each configuration, these tables give the number of converging instances\footnote{We must notice that some occurrences may not lead to a solution in the estimation procedure because $\Psi_{T_1,T_2}^n$ and $\Psi_{T_1,T_2,T_3}^n$, defined in Section~\ref{subsection:estimationTheta}, can sometimes take values outside the supports of $\psi_{T_1,T_2}^{-1}$ and $\psi_{T_1,T_2,T_3}^{-1}$.} of the estimator and their average, and the empirical quantile interval at $95\%$ (issued from taking the quantiles of the sample of estimated values). 
We observe that the estimators perform quite well: except on three lines in Table~\ref{tab:simulation4}, the true value of $\vartheta$ is always in the quantile interval.
Finally, we empirically observe that adding new maturities does not improve the quality of estimation in all configurations. For instance, increasing the number of maturities may increase or decrease the length of the quantile interval, and it may shift it away from the true value of $\vartheta$. Notice also that the one-step correction from $\widehat{\vartheta}_{2,n}$ to $\widetilde{\vartheta}_{2,n}$ never led to very different values.

\begin{table}[h]
        \centering
        \begin{tabular}{ccccc}
                \toprule
                Processes & Estimator & Instances that converged & Average & Quantile interval \\
                \midrule
                2 & $\widehat{\vartheta}_{2,n}$ & 100,000 & 1.4216 & [1.2697,1.6048] \\
                2 & $\widetilde{\vartheta}_{2,n}$ & 100,000 & 1.4217 & [1.2697,1.6048] \\
                3 & $\widehat{\vartheta}_{3,n}$ & 99,962 & 1.3799 & [0.77864,1.9250] \\
                4 & $\widehat{\vartheta}_{4,n}$ & 100,000 & 1.3840 & [1.0752,1.7646] \\
                5 & $\widehat{\vartheta}_{5,n}$ & 100,000 & 1.3807 & [1.1274,1.6864] \\
                6 & $\widehat{\vartheta}_{6,n}$ & 100,000 & 1.3849 & [1.0989,1.7644] \\
                \bottomrule
        \end{tabular}
        \caption{Results of the estimation on simulated data with $\vartheta=1.4$ yr$^{-1}$.}
        \label{tab:simulation1}
\end{table}

\begin{table}[h]
        \centering
        \begin{tabular}{ccccc}
                \toprule
                Processes & Estimator & Instances that converged & Average & Quantile interval \\
                \midrule
                2 & $\widehat{\vartheta}_{2,n}$ & 99,953 & 10.507 & [7.2997,16.500] \\
                2 & $\widetilde{\vartheta}_{2,n}$ & 99,953 & 10.507 & [7.2992,16.502] \\
                3 & $\widehat{\vartheta}_{3,n}$ & 100,000 & 9.9498 & [9.4916,10.258] \\
                4 & $\widehat{\vartheta}_{4,n}$ & 100,000 & 9.9424 & [9.6307,10.195] \\
                5 & $\widehat{\vartheta}_{5,n}$ & 100,000 & 9.9388 & [9.6538,10.180] \\
                6 & $\widehat{\vartheta}_{6,n}$ & 100,000 & 9.9511 & [9.6331,10.242] \\
                \bottomrule
        \end{tabular}
        \caption{Results of the estimation on simulated data with $\vartheta=10$ yr$^{-1}$}
        \label{tab:simulation2}
\end{table}

\begin{table}[h]
        \centering
        \begin{tabular}{ccccc}
                \toprule
                Processes & Estimator & Instances that converged & Average & Quantile interval \\
                \midrule
                2 & $\widehat{\vartheta}_{2,n}$ & 55,248 & 24.747 & [10.215,56.650] \\
                2 & $\widetilde{\vartheta}_{2,n}$ & 55,248 & 24.716 & [10.210,56.663] \\
                3 & $\widehat{\vartheta}_{3,n}$ & 100,000 & 33.904 & [22.598,40.060] \\
                4 & $\widehat{\vartheta}_{4,n}$ & 100,000 & 32.162 & [22.204,39.689] \\
                5 & $\widehat{\vartheta}_{5,n}$ & 100,000 & 31.075 & [22.046,38.832] \\
                6 & $\widehat{\vartheta}_{6,n}$ & 100,000 & 33.901 & [26.134,39.320] \\
                \bottomrule
        \end{tabular}
        \caption{Results of the estimation on simulated data with $\vartheta=40$ yr$^{-1}$.}
        \label{tab:simulation4}
\end{table}

\begin{figure}[h]
        \centering
        \includegraphics[width=13cm]{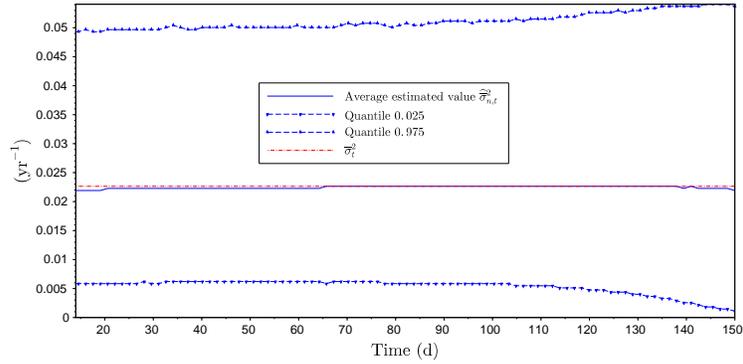}
        \caption{Quantiles for the square of the long-term volatility, with 2 processes, $\vartheta=10$ yr$^{-1}$ and deterministic constant volatilities}
        \label{fig:barresErreurNpSl2Processus}
\end{figure}

\begin{figure}[h]
        \centering
        \subfloat[Original plot]{\includegraphics[width=0.4\columnwidth]{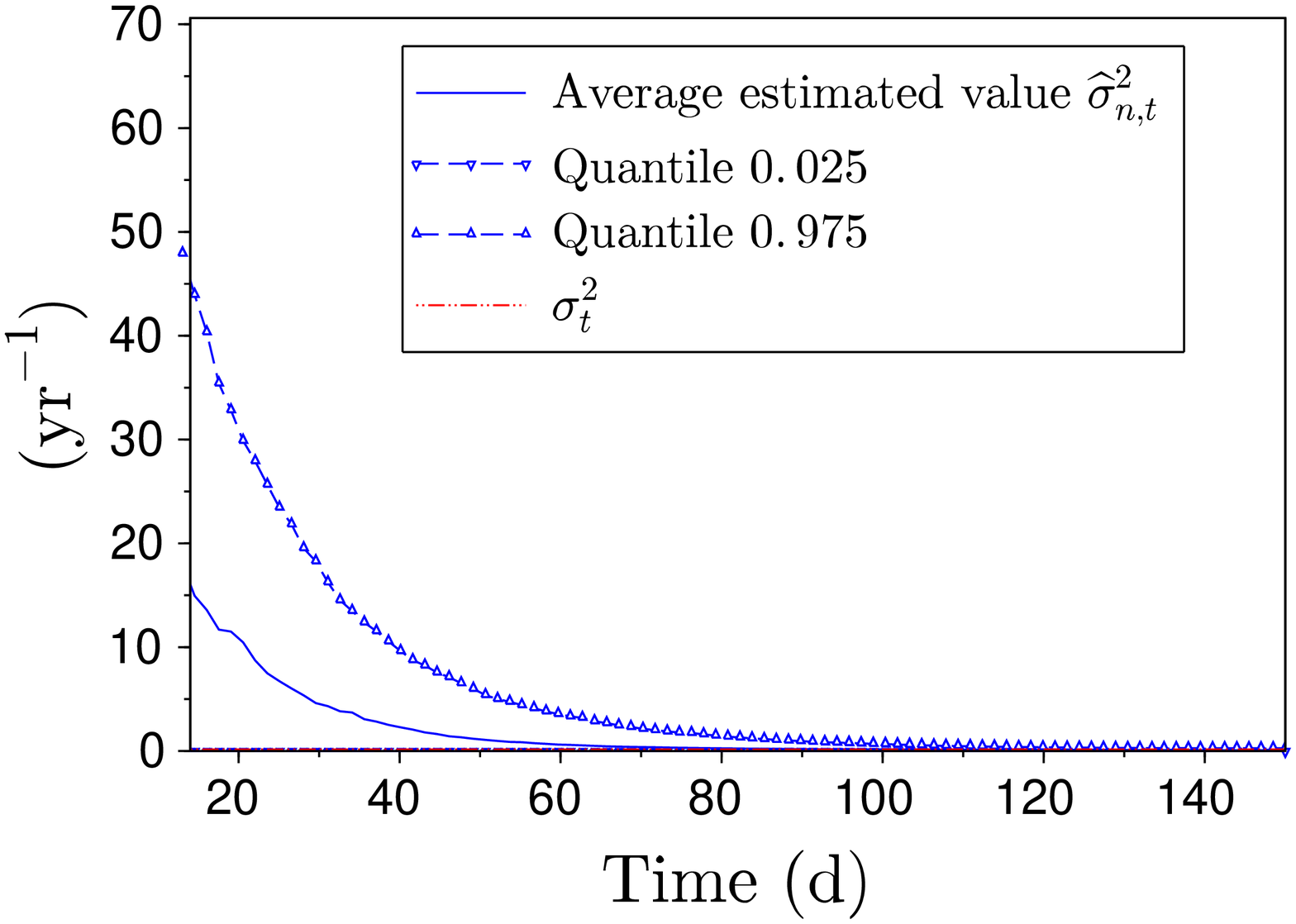}\label{fig:barresErreurNpSc2Processus}}
        \subfloat[Same plot with a ceiling of $1$~yr$^{-1}$ on the y-axis]{\includegraphics[width=0.4\columnwidth]{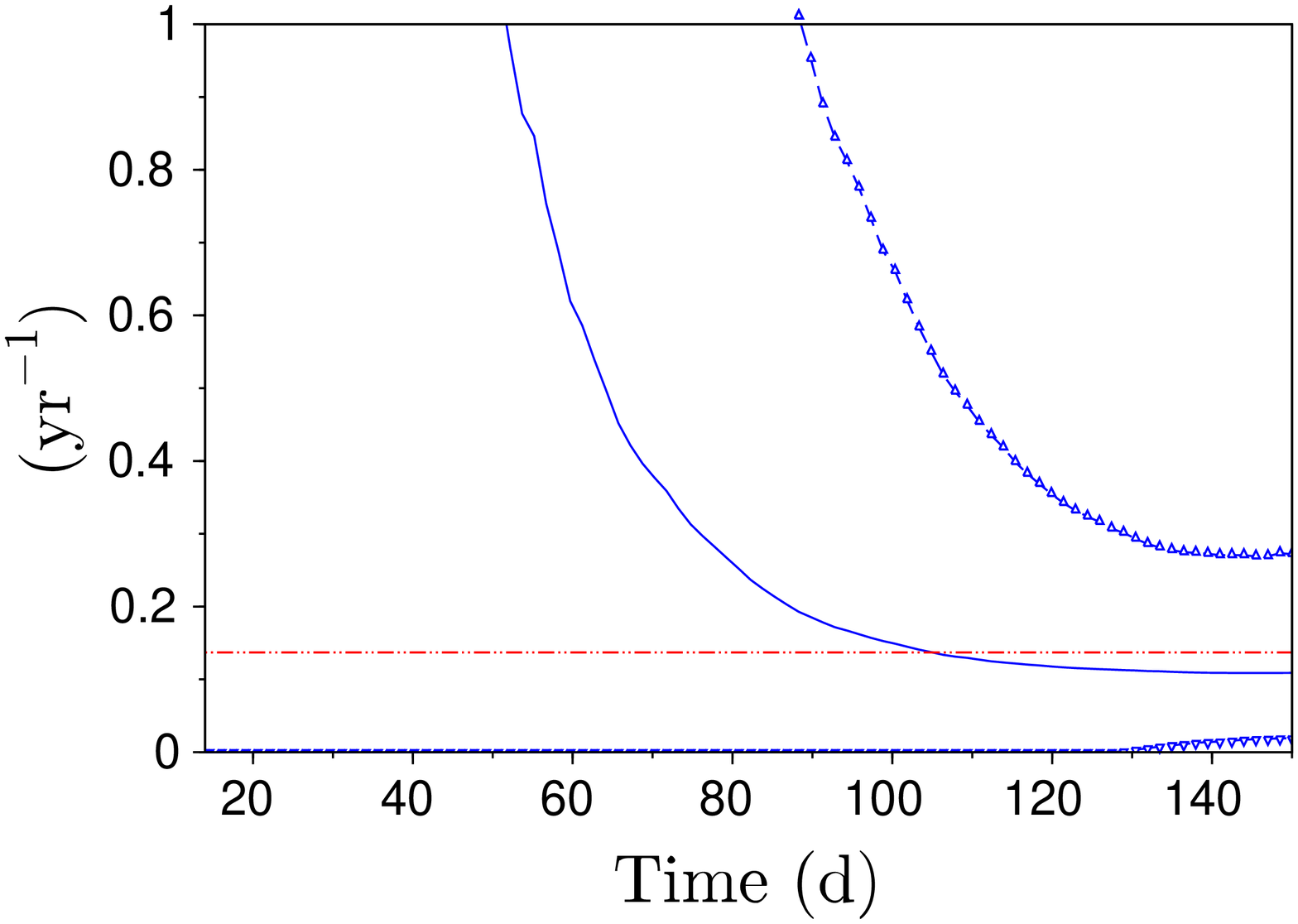}\label{fig:barresErreurNpSc2ProcessusZoom1}}\\
        \caption{Quantiles for the square of the short-term volatility, with 2 processes, $\vartheta=10$ yr$^{-1}$ and deterministic constant volatilities}
\end{figure}

Concerning the estimation results of the volatility processes $\sigmaS{t}^2$ and $\sigmaL{t}^2$, we use the causal kernel $K(x)=\indiq_{(0,1](x)}$, and the bandwidth $h_n$ for the two volatility functions is selected by cross validation. We also set $\varpi_n=3.65\cdot10^{-2}$~yr$^{-1}$. In the following we show the estimators $\widehat{\sigmaS{}}_n^2$ and $\widehat{\sigmaL{}}_n^2$ for the configuration where 2 processes are simulated on a period of 5 months (approximately 150 days), which means $T=T_1=150$ and $T_2=181$ days, with $n=100$ dates and $\vartheta = 10$~yr$^{-1}$. First we keep the specification $\drift{t}{j}=3.65\cdot 10^{-1}(\log(30)-X_t^j)$ for the drift process, but we use the constant volatility processes of \cite{bib:kiesel2009}, that is $\sigmaS{}=0.37$~yr$^{-1/2}$ and $\sigmaL{}=0.15$~yr$^{-1/2}$. A deterministic specification allows us to compare the curve of point estimates with the deterministic function that was used to simulate the processes. For this exercise with constant volatility functions, the cross validation criterion does not have a minimum, as taking the longest possible bandwidth is the best way to estimate a constant volatility. Therefore we use a bandwidth of $h_n=14$~d, which is the value that is selected when we implement the cross validation on real data, see below.
Remember that the nonparametric estimation result, Proposition~\ref{thm:estimationNPunePeriode}, gives convergence uniformly on $[h_n,T]$. Therefore we expect that the fit is not good for values ot $t$ being less than $h_n$, and this is why the lowest value on the x-axis is 14 days.
We perform simulation and estimation 100,000 times, and then take the average and the quantiles of the 100,000 curves (that is, at each point $t$ of the discretisation grid, we take the average and the quantiles at 2.5\% and 97.5\% of the 100,000 occurrences of $\widehat{\sigmaS{}}_{n,t}^2$ and $\widehat{\sigmaL{}}_{n,t}^2$). Figures~\ref{fig:barresErreurNpSl2Processus}  and \ref{fig:barresErreurNpSc2Processus} give the averages of $\widehat{\sigmaL{}}_{n,t}^2$ and $\widehat{\sigmaS{}}_{n,t}^2$, respectively, together with the true (constant) functions $\sigmaS{t}^2$ and $\sigmaL{t}^2$. The average of $\widehat{\sigmaS{}}_{n,t}^2$ is truncated, 0.4~\% of the values at each date being removed because of extreme values. The quantile intervals are plotted as well. They show a good estimation of the volatility function $\sigmaL{t}^2$. However, we can observe in Figure~\ref{fig:barresErreurNpSc2Processus} a bad performance of estimation of $\sigmaS{t}^2$, especially for large values of $T-t$, even when $t>h_n$. This reveals the ill-posedness of the problem, already pointed out in Section~\ref{subsection:estimationSigma}, due to the presence of the exponential terms in the denominator of the estimator of $\sigmaS{t}^2$: the latter can take very high values as soon as $\vartheta$ appears to be largely overestimated, which happens in a few simulations. Therefore, when high values appear, the estimation of $\sigmaS{t}$ should reasonably be taken into account only for small times to maturity $T-t$, where the estimation procedure seems to work well. Indeed, Figure~\ref{fig:barresErreurNpSc2ProcessusZoom1} shows that the quality of the fit of $\sigmaS{t}$ is reasonably good at the very right of the time range.
We put emphasis on the quantile intervals, that may be observed on all the plots. They show that the plot of the average exhibits a strange behaviour only because of some extreme values.

\par \medskip

\begin{figure}[h]
        \centering
        \includegraphics[width=13cm]{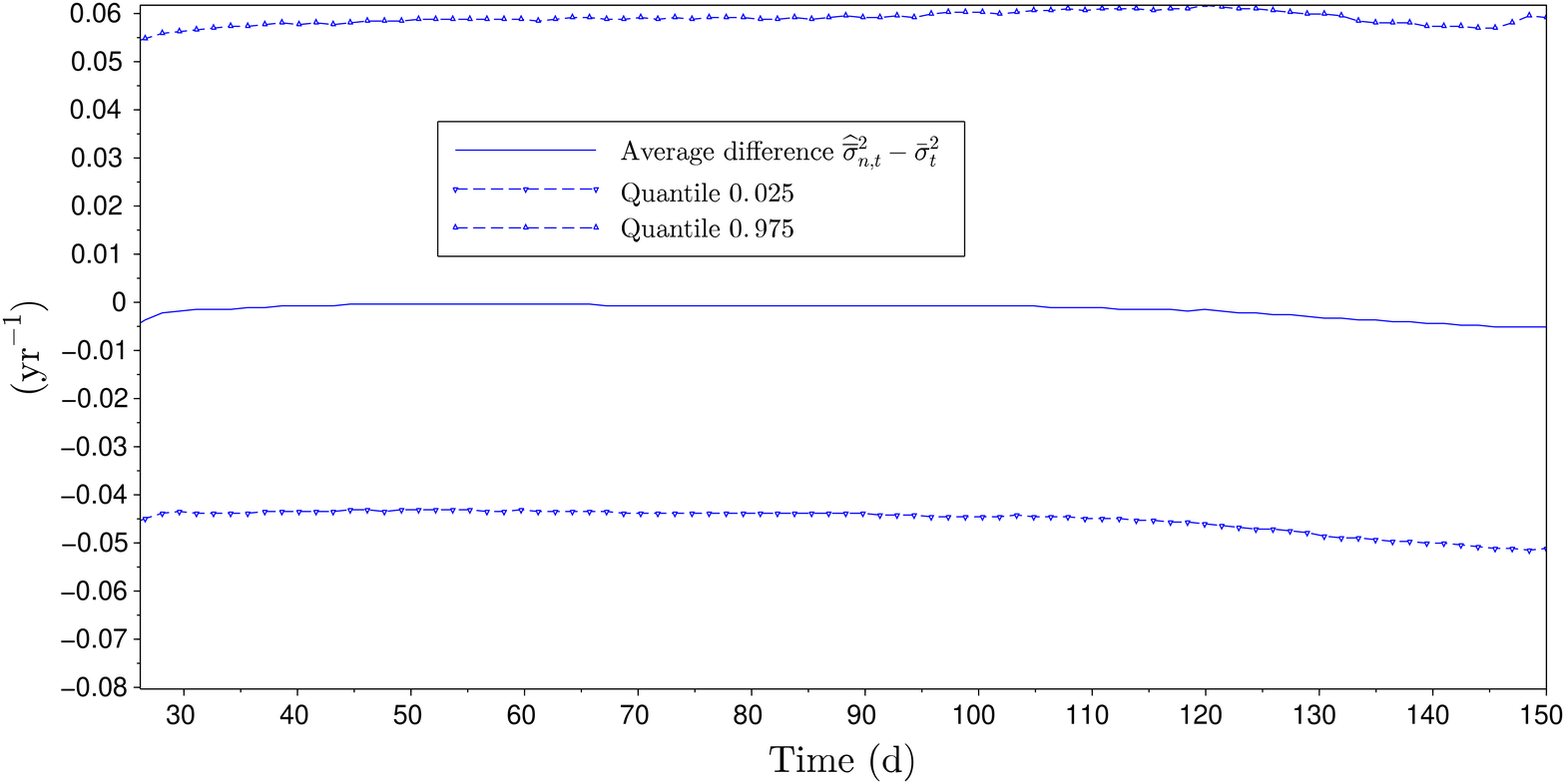}
        \caption{Estimation results of the long-term volatility process, in terms of $(\widehat{\overline{\sigma}}^2_{n,t}-\overline{\sigma}^2_t)$ with 2 processes, $\vartheta=10$ yr$^{-1}$ and the CIR-like specification for volatility processes}
        \label{fig:barresErreurNpSl2ProcessusVolSto}
\end{figure}

\begin{figure}[h]
        \centering
        \includegraphics[width=13cm]{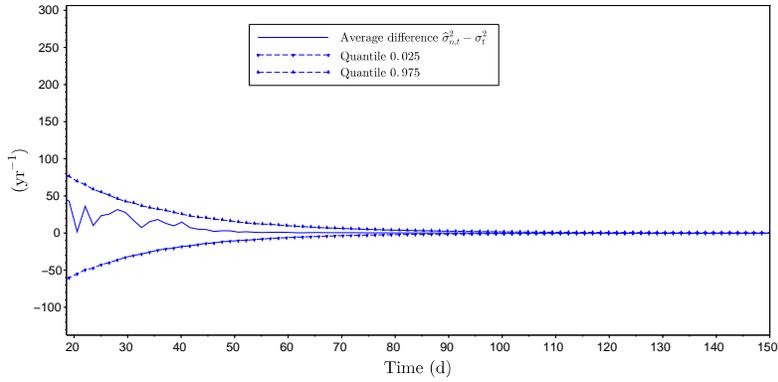}
        \caption{Estimation results of the short-term volatility process, in terms of $(\widehat{\sigma}^2_{n,t}-\sigma^2_t)$ with 2 processes, $\vartheta=10$ yr$^{-1}$ and the CIR-like specification for volatility processes}
        \label{fig:barresErreurNpSc2ProcessusVolSto}
\end{figure}

Now, we are back to the specification $\sigmaS{t}=0.37 \Sigma_t^d, \sigmaL{t}=0.15 \Sigma_t^d$. As the volatility processes depend on the path of $X$, we cannot compare visually the real volatility and its point estimators. This is why we plot the differences $\widehat{\sigmaS{}}_{n,t}^2-\sigmaS{t}^2$ and $\widehat{\sigmaL{}}_{n,t}^2-\sigmaL{t}^2$ instead. Moreover, we let $\widehat{\sigmaS{}}_{n,t}^2$ and $\widehat{\sigmaL{}}_{n,t}^2$ take negative values to emphasize symmetry in the differences. We plot the average and the quantile curves of the 100,000 estimators for the two volatility processes in Figures~\ref{fig:barresErreurNpSl2ProcessusVolSto} and \ref{fig:barresErreurNpSc2ProcessusVolSto}. The behaviours of the series of point estimators are very similar to the ones we described while considering deterministic volatility functions. The bandwidths for nonparametric estimation are chosen using cross validation on each simulation. The average bandwidths are 18.6 days for $\sigmaS{}$ and 26.2 days for $\sigmaL{}$.
Because the estimator $\widehat{\sigmaS{}}_{n,t}^2$ can take very high values, as we have noticed above, the average that is plotted in Figure~\ref{fig:barresErreurNpSc2ProcessusVolSto} is computed using the 100,000 point estimates available at each plotting date, except the 100 lowest and the 100 highest ones. The quantiles are still plotted using all the available point estimates at each plotting date.

\subsection{Results on real data from the French electricity market} \label{sec_french_electricity_market}

The data used for estimation are the 6 available month-ahead forward contracts on the French market (\url{www.eex.com}) from December 6\textsuperscript{th}, 2001 to November 30\textsuperscript{th}, 2018. On this history, we get 204 periods of 1 month where 6 processes (the 6 month-ahead contracts) are jointly observed, whereas we get 200 periods of 5 months where 2 processes (the 1 month-ahead and the 2 month-ahead contracts) are jointly observed. These numbers of periods are given in Table~\ref{tab:estimateursSurDonneesReellesFrance} for all the configurations described in Section~\ref{sec_res_simulation_data}. In the same column, Table~\ref{tab:estimateursSurDonneesReellesFrance} also precises the number of periods on which the estimator converges. And the same table gives the estimation results of $\vartheta$ for all the possible configurations, with the average value and the standard deviation of the estimators.
\begin{table}[h]
        \centering
        \begin{tabular}{cccc}
                \toprule
                Estimator & Per. with convergence/ Number of per. & Average & Standard deviation \\
                \midrule
                $\widehat{\vartheta}_{2,n}$ & 71/200 & 24.966 & 10.734 \\
                $\widetilde{\vartheta}_{2,n}$ & 71/200 & 25.049 & 10.743 \\
                $\widehat{\vartheta}_{3,n}$ & 142/201 & 4.488 & 3.462 \\
                $\widehat{\vartheta}_{4,n}$ & 149/202 & 3.387 & 2.732 \\
                $\widehat{\vartheta}_{5,n}$ & 149/203 & 3.072 & 2.800 \\
                $\widehat{\vartheta}_{6,n}$ & 132/204 & 3.514 & 2.946 \\
                \bottomrule
        \end{tabular}
        \caption{Estimators of $\vartheta$ on real data in France (unit: yr$^{-1}$).}
        \label{tab:estimateursSurDonneesReellesFrance}
\end{table}
Contrary to the results on simulated data, the values of the estimators are different from one configuration to another. More precisely, the estimators from 2 processes are higher (of a factor between 5 and 8) than the ones from 3 to 6 processes. This can be explained by two different causes. First, we have proved the convergence of the estimators from 3 to 6 processes only if the drift is zero: this was stated in Theorem~\ref{thm:loiAsymptotiqueHatAMoinsA}. Second, these differences may be due to the presence of errors linked to measurement or to the model. When we compare the estimators from more than 3 processes, we can observe that, like in the case of simulated data, adding processes does not significantly improve the estimation results. 

We end this session by performing estimation on the last part of the dataset, which runs from July~1\textsuperscript{st}, 2018 to November~30\textsuperscript{th}, 2018. The value taken by $\widehat{\vartheta}_{2,n}$ is 32.781 yr$^{-1}$, while $\widetilde{\vartheta}_{2,n}$ has value 32.779 yr$^{-1}$. We use Theorems~\ref{thm:loiAsymptotiqueHatAMoinsA} and \ref{thm:estimationEfficace} to provide feasible central limit theorems, and thus we get confidence intervals using the nonparametric estimates of the volatility processes. The confidence intervals are $[28.814,36.748]$ for $\widehat{\vartheta}_{2,n}$ and $[29.918,35.640]$ for $\widetilde{\vartheta}_{2,n}$ (values in yr$^{-1}$). As expected because of efficiency, the latter is narrower than the former.\par
\begin{figure}[h]
        \centering
        \subfloat[Plot of $\widehat{\overline{\sigma}}_{n,t}^2$ obtained with the estimator $\widehat{\vartheta}_{2,n}$]{\includegraphics[width=0.4\columnwidth]{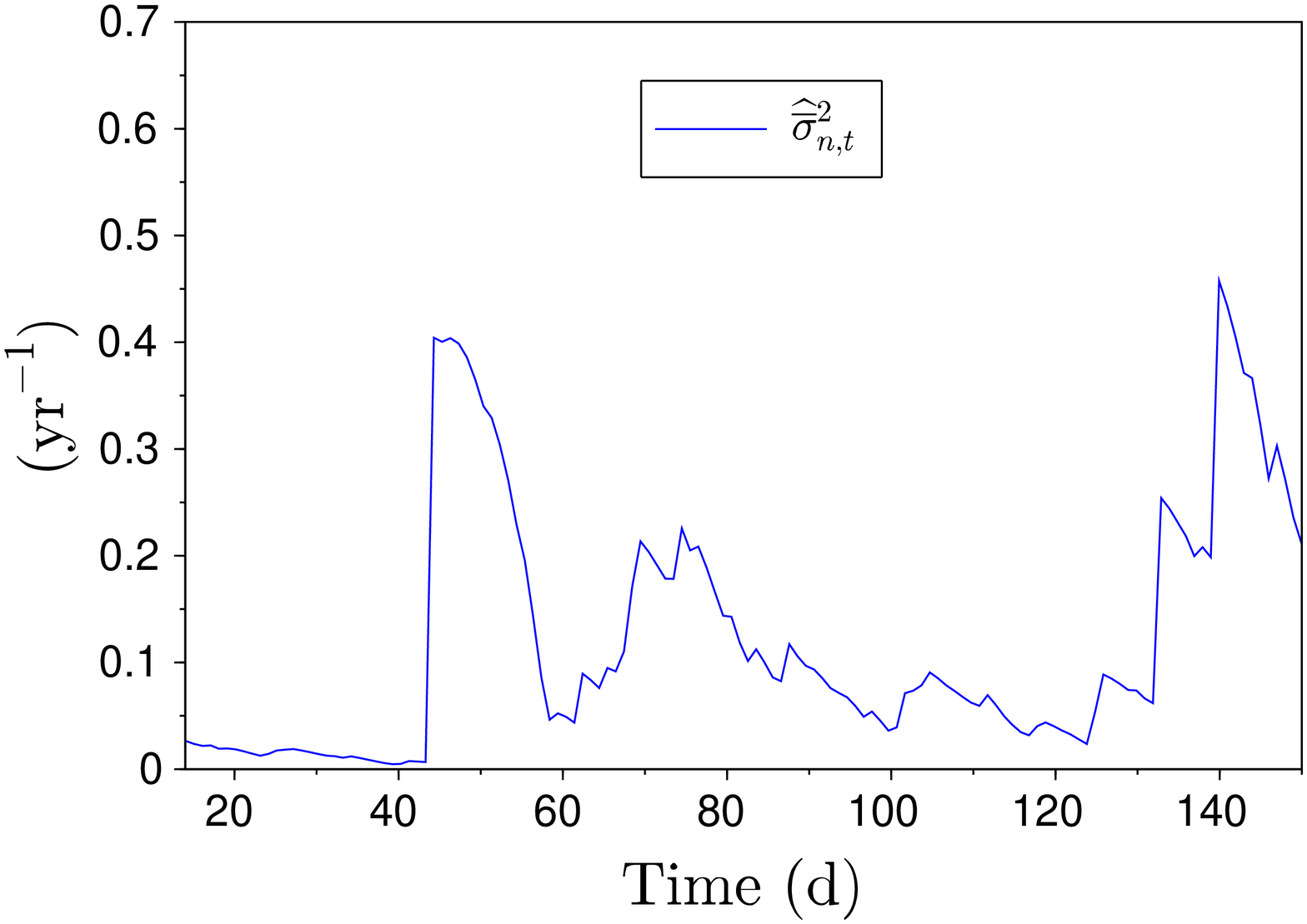}}
        \subfloat[Plot of $\widehat{\overline{\sigma}}_{n,t}^2$ obtained with the estimator $\widehat{\vartheta}_{3,n}$]{\includegraphics[width=0.4\columnwidth]{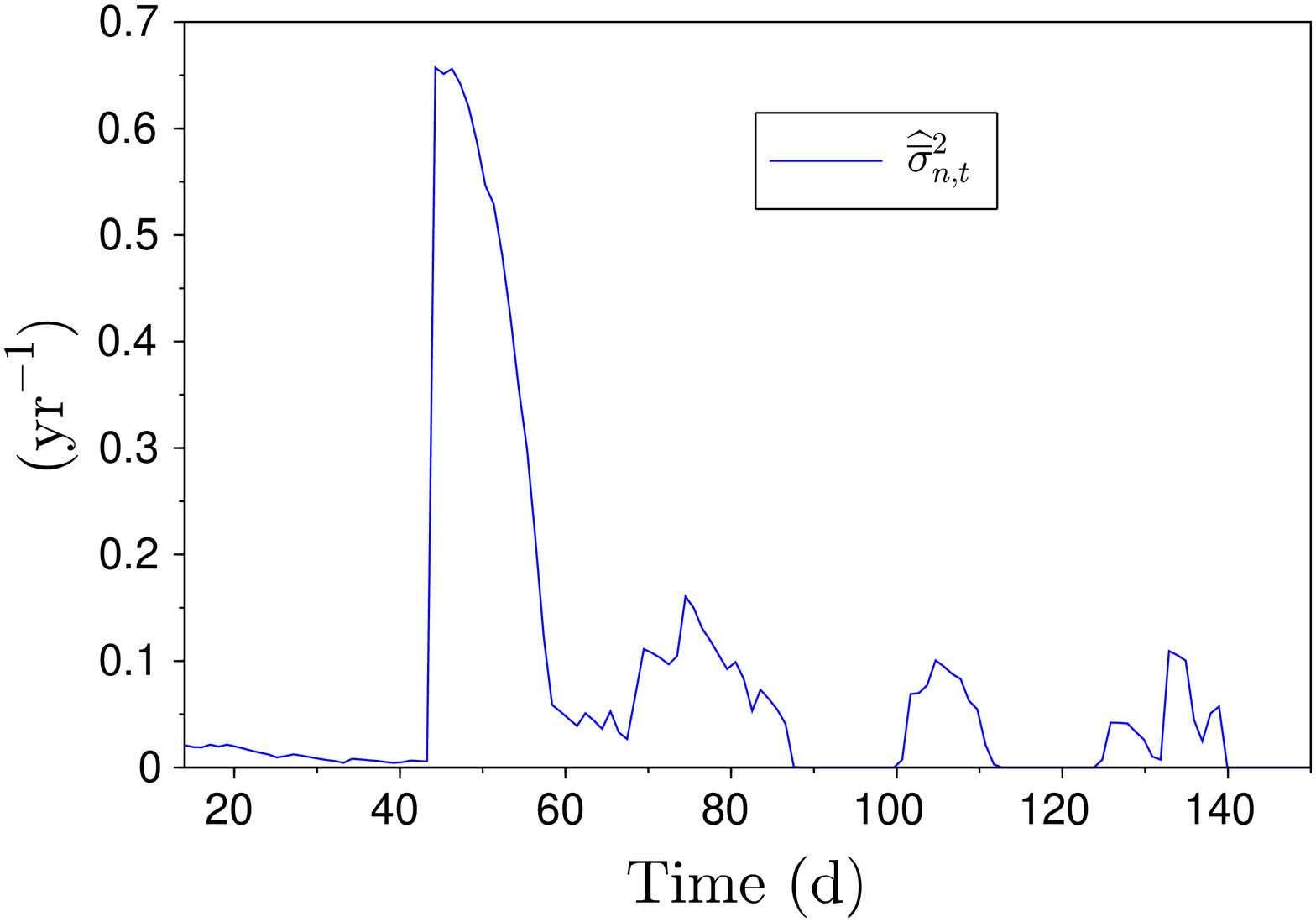}}\\
        \subfloat[Plot of $\widehat{\sigma}_{n,t}^2$ obtained with the estimator $\widehat{\vartheta}_{2,n}$]{\includegraphics[width=0.4\columnwidth]{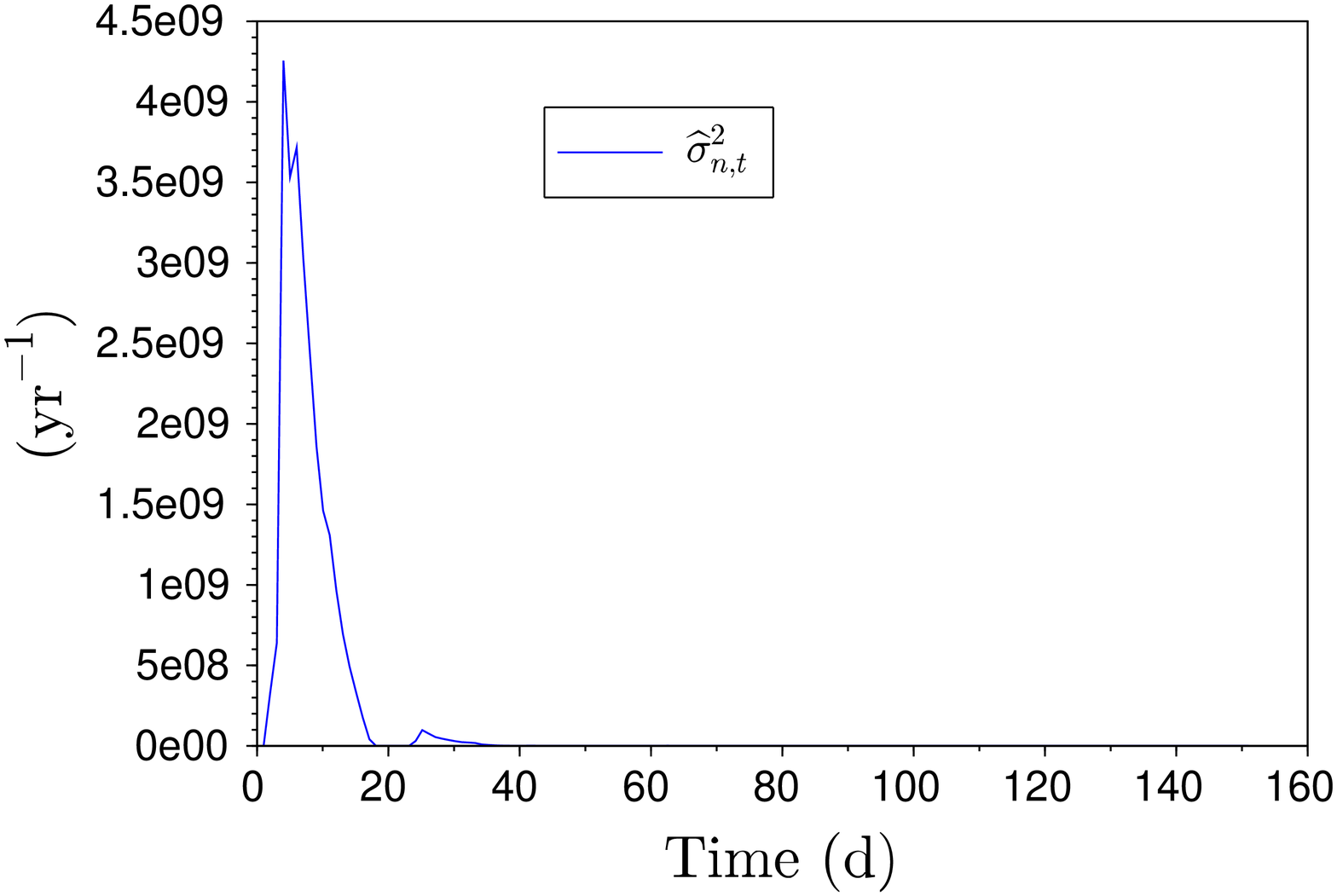}}
        \subfloat[Plot of $\widehat{\sigma}_{n,t}^2$ obtained with the estimator $\widehat{\vartheta}_{3,n}$]{\includegraphics[width=0.4\columnwidth]{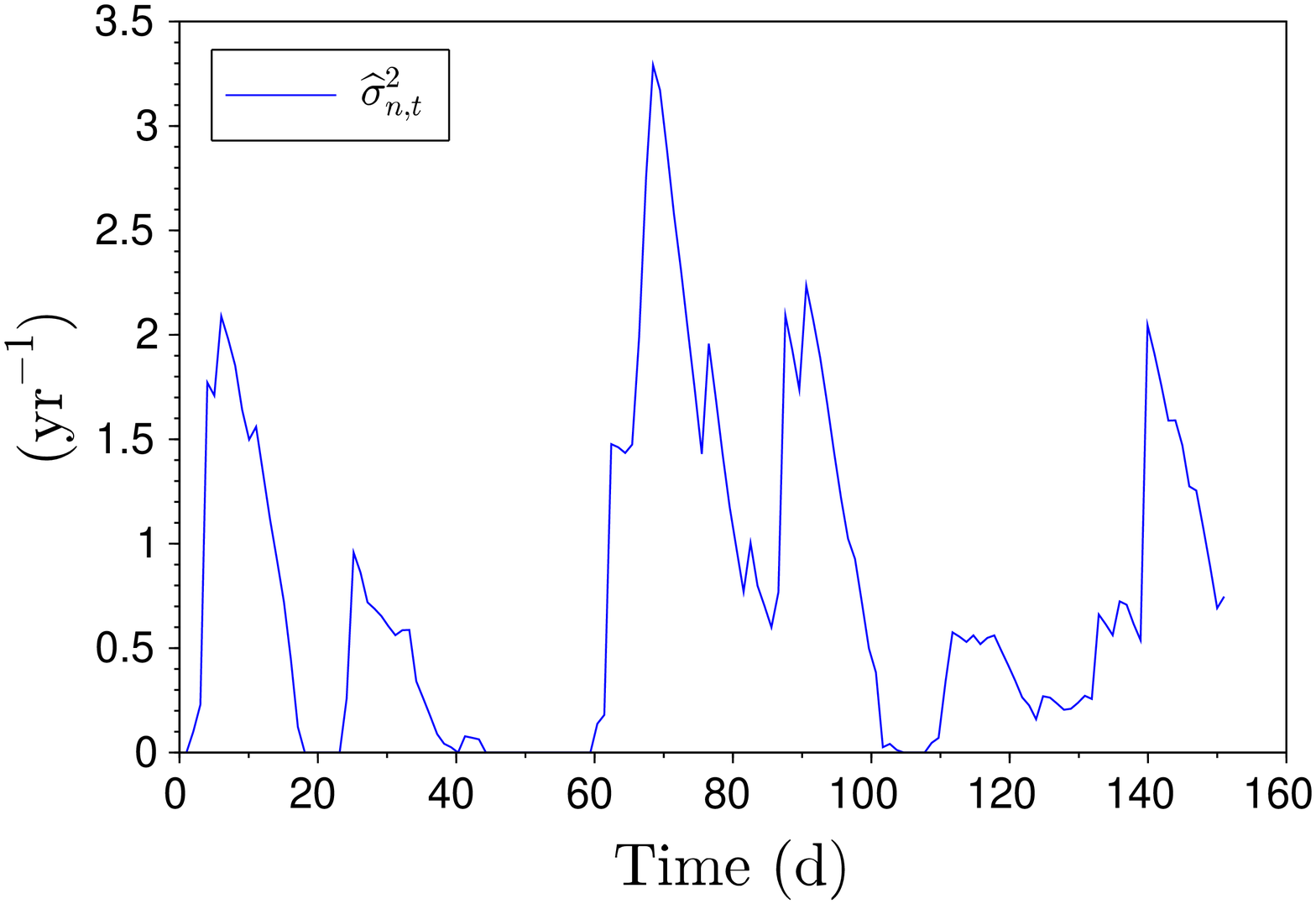}}\\
        \subfloat[Plot of $e^{-2\widehat{\vartheta}_{2,n}(T_1-t)}\widehat{\sigma}_{n,t}^2$ obtained with the estimator $\widehat{\vartheta}_{2,n}$]{\includegraphics[width=0.4\columnwidth]{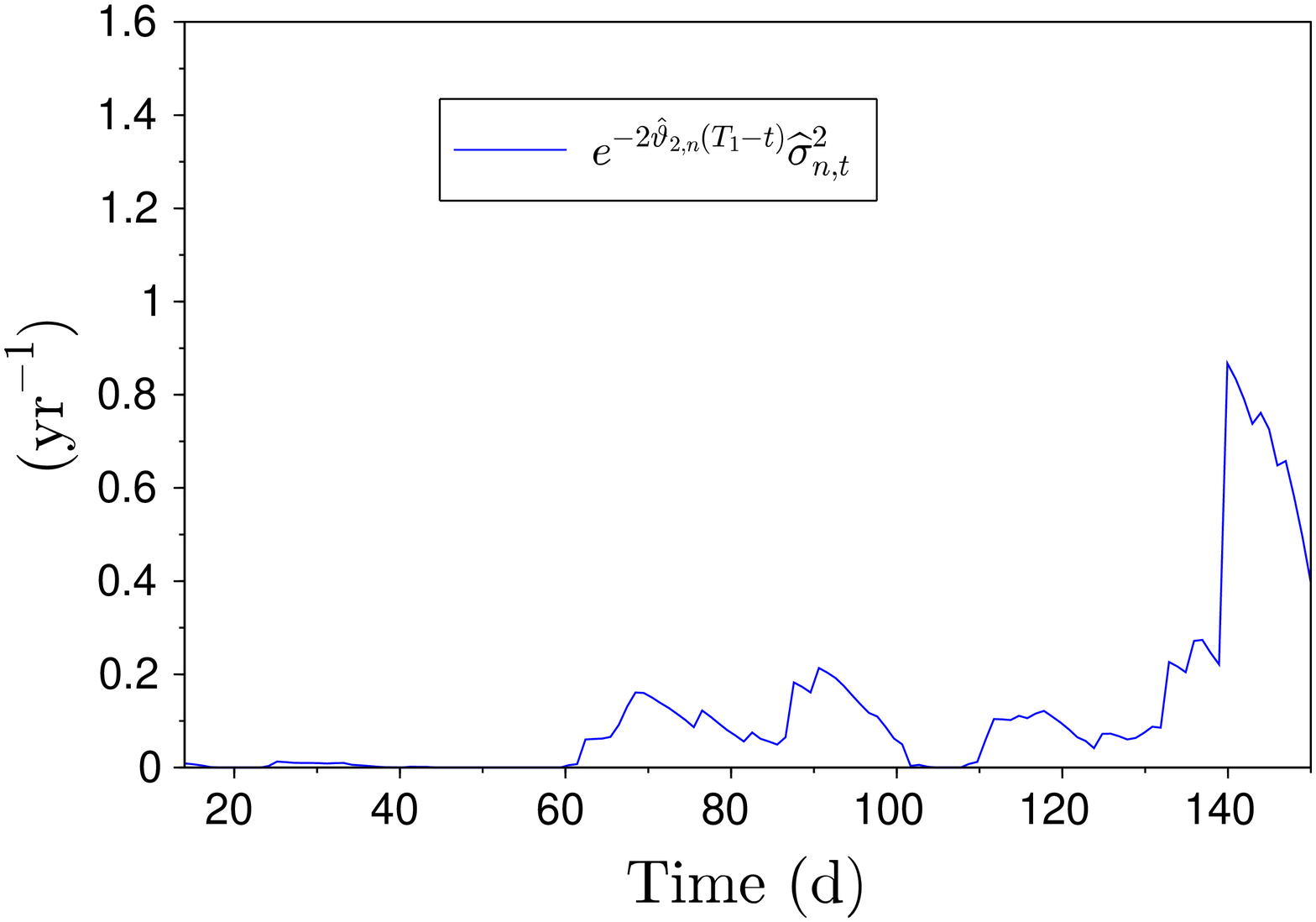}}
        \subfloat[Plot of $e^{-2\widehat{\vartheta}_{3,n}(T_1-t)}\widehat{\sigma}_{n,t}^2$ obtained with the estimator $\widehat{\vartheta}_{3,n}$]{\includegraphics[width=0.4\columnwidth]{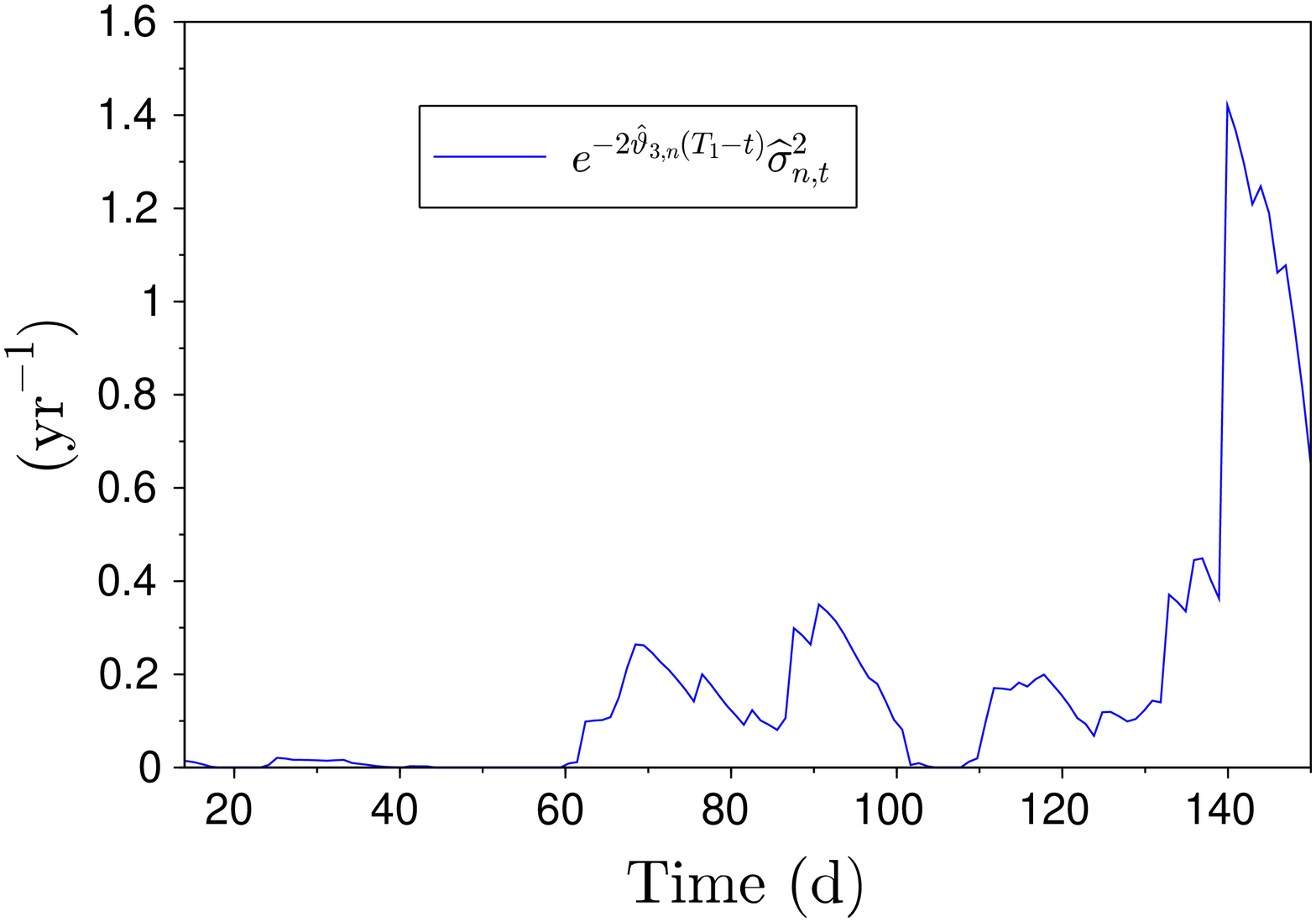}}
        \caption{Plots of estimated paths of volatility processes. Left: using $\widehat{\vartheta}_{2,n}$ as a preliminary estimator of $\vartheta$. Right: using $\widehat{\vartheta}_{3,n}$.}
        \label{fig:realdata_volatility}
\end{figure}

We plot the estimates of the two volatility processes in Figure~\ref{fig:realdata_volatility}. The bandwidths, chosen using cross validation, are equal to 14 days. These results clearly show non constant volatility processes. In both cases, the estimates of the volatility process $\widehat{\overline{\sigma}}_{n,t}$ evolve in a same order of magnitude, although we can observe higher values of the estimate $\widehat{\overline{\sigma}}_{n,t}$ from $\widehat{\vartheta}_{2,n}$, in the last 70 days. The difference between $\widehat{\vartheta}_{2,n}$ and $\widehat{\vartheta}_{3,n}$ has a higher impact on the estimation of $\widehat{\sigma}_{n,t}^2$. And, the high value of $\widehat{\vartheta}_{2,n}$ leads to huge values of $\widehat{\sigma}_{n,t}^2$ especially when $T_1-t$ is large, due to the exponential term in the short term volatility, as explained in Section \ref{sec_res_simulation_data}. However, the estimates of the global volatility function $e^{-2\vartheta(T_1-t)}\sigma^2$ are close in both cases, although taking lower values in the case of $\widehat{\vartheta}_{2,n}$. This experiment reveals, like on simulated data, the ill-posedness of the problem of estimation of the volatility paths, especially for $\sigmaS{}$ when $\widehat{\vartheta}_{2,n}$ is large. 

\section{Proofs} \label{section:preuves}

\subsection{Preliminaries: localisation} \label{localization}

With no loss of generality, we may (and will) assume that the processes $b^j$, $\sigmaS{}$ and $\sigmaL{}$ are bounded, relying on a so-called localisation argument. For an integer $p \geq 1$, introduce the stopping time $\tau_p = \inf\{t \in [0,T], \min(b^j_t,\sigma_t, \overline{\sigma}_t) > p\}$. Replacing $X_t^j$ by $X_{t \wedge \tau_p}$, we have bounded processes $b^j$, $\sigmaS{}$ and $\sigmaL{}$. Moreover since these processes are at least locally bounded, we have $\mathbb P(\tau_p > T) \rightarrow 1$ as $p \rightarrow \infty$. 
We refer to Section 3.6.3 in \cite{bib:statisticsAndHFDataInStatisticalMethodsForSDE} for details. 

\subsection{Proof of Theorem~\ref{thm:loiAsymptotiqueHatAMoinsA}} \label{subsection:preuveToutPremierTheoreme}

\subsubsection*{Proof of Theorem \ref{thm:loiAsymptotiqueHatAMoinsA} (1)} \label{proof of thm 1}

\noindent {\bf Step 1.} We first consider the case $\drift{}{j}=0$ for $j=1,2$. For notational simplicity, we set $\diffE{\ell}{k}{\vartheta} = e^{-\vartheta T_k}-e^{-\vartheta T_\ell}$. Let us define
\[
        \zeta_i^n  = (\Delta_i^n X^2)^2-(\Delta_i^n X^1)^2
\]
and
\[
        \xi_i^n  = \big(\Delta_i^n X^2-\Delta_i^n X^1\big)^2.
\]
Clearly
\begin{align*}
  & \Big(\int_{(i-1)\Delta_n}^{i\Delta_n}e^{-\vartheta(T_2-t)}\sigmaS{t}d\Ws{t}+\int_{(i-1)\Delta_n}^{i\Delta_n}\sigmaL{t}d\Wl{t}\Big)^2 
          -\Big(\int_{(i-1)\Delta_n}^{i\Delta_n}e^{-\vartheta(T_1-t)}\sigmaS{t}d\Ws{t}+\int_{(i-1)\Delta_n}^{i\Delta_n}\sigmaL{t}d\Wl{t}\Big)^2 \\
        = & (e^{-2\vartheta T_2}-e^{-2\vartheta T_1}) \Big( \int_{(i-1)\Delta_n}^{i\Delta_n}e^{\vartheta t}\sigmaS{t}d\Ws{t}\Big)^2 
        +2\diffE{1}{2}{\vartheta} \int_{(i-1)\Delta_n}^{i\Delta_n}e^{\vartheta t}\sigmaS{t}d\Ws{t}\int_{(i-1)\Delta_n}^{i\Delta_n}\sigmaL{t}d\Wl{t},
\end{align*}
therefore, setting $\chi_i^n  = 2\diffE{1}{2}{\vartheta} \int_{(i-1)\Delta_n}^{i\Delta_n}e^{\vartheta t}\sigmaS{t}d\Ws{t}\int_{(i-1)\Delta_n}^{i\Delta_n}\sigmaL{t}d\Wl{t}$, we obtain the following representation
\begin{equation} \label{key rep}
\zeta_i^n =  \frac{1}{\psi_{T_1,T_2}(\vartheta)}\xi_i^n + \chi_i^n.
\end{equation}
By standard convergence of the quadratic variation (see for instance Section~2.1.5 in \cite{bib:econometricsOfHFDataInStatisticalMethodsForSDE}),
\[
        \sum_{i = 1}^n \xi_i^n \rightarrow \diffE{1}{2}{\vartheta}^2\int_0^T e^{2\vartheta t}\sigmaS{t}^2d{t}
\]
in probability. Note that the limit is almost surely positive by Assumption~\ref{asm:regulariteFonctions}. Also, since $\Ws{}$ and $\Wl{}$ are independent, and since $\sigmaS{t}^2 \leq M$ and $\sigmaL{t}^2 \leq M$ for some constant $M>0$ by localization, we have that 
\begin{align*}
& \E\Big[\Big(\sum_{i = 1}^n \int_{(i-1)\Delta_n}^{i\Delta_n}e^{\vartheta t}\sigmaS{t}d\Ws{t}\int_{(i-1)\Delta_n}^{i\Delta_n}\sigmaL{t}d\Wl{t}\Big)^2\Big]  \\
= & \sum_{i = 1}^n\E\big[\big(\int_{(i-1)\Delta_n}^{i\Delta_n}e^{\vartheta t}\sigmaS{t}d\Ws{t}\big)^2\big]\E\big[\big(\int_{(i-1)\Delta_n}^{i\Delta_n}\sigmaL{t}d\Wl{t}\big)^2\big] \leq \Delta_n e^{2\vartheta T}M^2 \rightarrow 0.
\end{align*}
Therefore $\sum_{i = 1}^n \zeta_i^n$ converges in probability as well, with the same limit as $\frac{1}{\psi_{T_1,T_2}(\vartheta)}\sum_{i = 1}^n\xi_i^n$. It follows that
\[
        \Psi_{T_1,T_2}^n  = \frac{\sum_{i = 1}^n \big(\Delta_i^n X^2-\Delta_i^n X^1 \big)^2}{\sum_{i = 1}^n (\Delta_i^n X^2)^2-(\Delta_i^n X^1)^2} = \frac{\sum_{i=1}^n \xi_i^n}{\sum_{i = 1}^n \zeta_i^n} \rightarrow  \psi_{T_1,T_2}(\vartheta)
\]
in probability. We derive the convergence
\[
        \psi_{T_1,T_2}(\widehat{\vartheta}_{2,n}) \rightarrow \psi_{T_1,T_2}(\vartheta)
\]
in probability on the event $\{\Psi_{T_1,T_2}^n \in (-1,0)\}$, hence the convergence $\widehat{\vartheta}_{2,n} \rightarrow \vartheta$ in probability as well since $\{\Psi_{T_1,T_2}^n \in (-1,0)\}$ has asymptotically probability $1$ and that $\vartheta \leadsto \psi_{T_1, T_2}(\vartheta)$ is invertible with continuous inverse.\par \bigskip

\noindent {\bf Step 2.} Using \eqref{key rep}, we readily obtain
\begin{align*}
\Delta_n^{-1/2}\big(\Psi_{T_1,T_2}^n - \psi_{T_1,T_2}(\vartheta)\big) & = \Delta_n^{-1/2}\Big(\frac{\sum_{i=1}^n \xi_i^n}{\sum_{i = 1}^n \zeta_i^n} - \psi_{T_1,T_2}(\vartheta)\Big) 
 = -\psi_{T_1,T_2}(\vartheta) \frac{\Delta_n^{-1/2}\sum_{i=1}^n \chi_i^n}{\sum_{i = 1}^n \zeta_i^n}.
\end{align*}
Consider next the sequence of 1-dimensional processes
\[
        \chi_n(t) = \Delta_n^{1/2}\sum_{i = 1}^{\lfloor t\Delta_n^{-1}\rfloor}f\big(\Delta_n^{-1/2}\Delta_i^n Y^1,\Delta_n^{-1/2}\Delta_i^n Y^2\big),
\]
where $Y_t = (Y_t^1,Y_t^2) = \big(\int_0^t e^{\vartheta s}\sigmaS{s} d\Ws{s}, \int_0^t \sigmaL{s}d\Wl{s}\big)$. By Theorem~3.21, p.~231 in \cite{bib:statisticsAndHFDataInStatisticalMethodsForSDE} applied to the martingale $Y$ with $f(x,y)=xy$ which has vanishing integral under the standard 2-dimensional-Gaussian measure, we have that the process $\chi_n(t)$ converges stably in law to a continuous process $\chi(t)$ defined on an extension of the original probability space and given by
\[
        \chi(t) = \int_0^t e^{\vartheta s} \sigmaS{s} \sigmaL{s} dW_s,
\]
where $W$ is a Brownian motion independent of $\cF$. Using successively $\Delta_n^{-1/2}\sum_{i = 1}^n \chi_i^n =  2\diffE{1}{2}{\vartheta} \chi_n(T)$, the fact that the convergence $\chi_n\rightarrow \chi$ holds stably in law and the convergence 
\[
        \sum_{i = 1}^n \zeta_i^n \rightarrow (e^{-2\vartheta T_2}-e^{-2\vartheta T_1})\int_0^T e^{2\vartheta t} \sigmaS{t}^2 dt,
\]
in probability, we derive 
\begin{align*}
-\Delta_n^{-1/2}\psi_{T_1,T_2}(\vartheta) \frac{\sum_{i=1}^n \chi_i^n}{\sum_{i = 1}^n \zeta_i^n}  \rightarrow  & -\psi_{T_1,T_2}(\vartheta) \frac{2(e^{-\vartheta T_2}-e^{-\vartheta T_1}\big) \chi(T)}{(e^{-2\vartheta T_2}-e^{-2\vartheta T_1})\int_0^T e^{2\vartheta t} \sigmaS{t}^2 dt} \\
=  & -\frac{2(e^{-\vartheta T_2}-e^{-\vartheta T_1})^3}{(e^{-2\vartheta T_2}-e^{-2\vartheta T_1})^2\int_0^T e^{2\vartheta t} \sigmaS{t}^2 dt}\chi(T)
\end{align*}
in distribution. Conditional on $\cF$, the limiting variable is centred Gaussian, with conditional variance $v_\vartheta(\sigmaS{}, \sigmaL{}) = 4\frac{(e^{-\vartheta T_2}-e^{-\vartheta T_1})^2}{(e^{-\vartheta T_2}+e^{-\vartheta T_1})^4} \frac{\int_0^T e^{2\vartheta t}\sigmaS{t}^2\sigmaL{t}^2dt}{(\int_0^T e^{2\vartheta t}\sigmaS{t}^2dt)^2}$.\par \bigskip

\noindent {\bf Step 3.} On the event $\{\Psi_{T_1,T_2}^n \in (-1,0)\}$, we have 
\[
        \Delta_n^{-1/2}\big(\widehat{\vartheta}_{2,n} -\vartheta\big) = \Delta_n^{-1/2}\big(\Psi_{T_1,T_2}^n - \psi_{T_1,T_2}(\vartheta)\big)\partial_\vartheta \psi_{T_1,T_2}^{-1}(Z_n)
\]
for some $Z_n$ that converges to $\psi_{T_1,T_2}(\vartheta)$ in probability by Step~1. The conclusion follows from 
\[
        \big(\partial_\vartheta \psi_{T_1,T_2}^{-1}(\psi_{T_1,T_2}(\vartheta))\big)^2v_\vartheta(\sigmaS{}, \sigmaL{}) = V_\vartheta(\sigmaS{}, \sigmaL{})
\]
together with the fact that $\{\Psi_{T_1,T_2}^n \in (-1,0)\}$ has asymptotically probability $1$.\par \bigskip

\noindent {\bf Step 4.} It remains to relax  the restriction $\drift{}{j}=0$. When $\drift{}{j}$ is non-zero, by localization again, we may assume it is bounded. Moreover, the volatility processes $(\sigma, \overline{\sigma})$ are bounded below by Assumption \ref{asm:regulariteFonctions}. Then, by Girsanov theorem, we apply a change of measure which is $\cF_T$-measurable. Since the convergence in distribution in Step~2 holds stably in law, we may work under this change of measure (see Section~2.4.4 in \cite{bib:econometricsOfHFDataInStatisticalMethodsForSDE} for a simple explanation)). Finally, relaxing the boundedness assumption on $\sigmaS{}, \sigmaL{}$ and $\drift{}{j}$ is standard, see Section~\ref{localization} above.

\subsubsection*{Proof of Theorem \ref{thm:loiAsymptotiqueHatAMoinsA} (2)}

\noindent {\bf Step 1.} We have
\[
        \Psi_{T_1,T_2,T_3}^n  = \frac{\sum_{i =1}^n \big(\Delta_i^n(X^3-X^2)\big)^2}{\sum_{i = 1}^n \big(\Delta_i^n(X^2-X^1)\big)^2}
\]
By standard convergence of the quadratic variation
\begin{align}
        & \sum_{i = 1}^n \big(\Delta_i^n(X^2-X^1)\big)^2 \rightarrow \diffE{1}{2}{\vartheta}^2\int_0^Te^{2\vartheta t}\sigmaS{t}^2 dt,  \label{quad bis}\\
        & \sum_{i = 1}^n \big(\Delta_i^n(X^3-X^2)\big)^2 \rightarrow \diffE{2}{3}{\vartheta}^2\int_0^T e^{2\vartheta t}\sigmaS{t}^2 dt \nonumber
\end{align}
in probability. Since $\psi_{T_1,T_2,T_3}(\vartheta) = \frac{\diffE{2}{3}{\vartheta}^2}{\diffE{1}{2}{\vartheta}^2}$, we derive 
$\psi_{T_1,T_2,T_3}\big(\widehat{\vartheta}_{n,3}\big) \rightarrow \psi_{T_1,T_2,T_3}(\vartheta)$
in probability on the event $\big\{\Psi_{T_1,T_2,T_3}\in \big(0, \big(\frac{T_3-T_2}{T_2-T_1}\big)^2\big)\big\}$ which has asymptotically probability $1$, hence the convergence $\widehat{\vartheta}_{n,3} \rightarrow \vartheta$ in probability.\par \bigskip

\noindent {\bf Step 2.} We further have
\[
        \Psi_{T_1,T_2,T_3}^n -\psi_{T_1,T_2,T_3}(\vartheta) = \frac{\sum_{i =1}^n \big(\Delta_i^n(X^3-X^2)\big)^2}{\sum_{i = 1}^n \big(\Delta_i^n(X^2-X^1)\big)^2}- \frac{\diffE{2}{3}{\vartheta}^2}{\diffE{1}{2}{\vartheta}^2} = \frac{\sum_{i =1}^n \eta_i^n}{\sum_{i = 1}^n \big(\Delta_i^n(X^2-X^1)\big)^2},
\]
with 
\[
        \eta_i^n=\big(\Delta_i^n(X^3-X^2)\big)^2-\frac{\diffE{2}{3}{\vartheta}^2}{\diffE{1}{2}{\vartheta}^2}\big(\Delta_i^n(X^2-X^1)\big)^2.
\]
Write $\overline{\Delta}_i^nf=\int_{(i-1)\Delta_n}^{i\Delta_n}f(t)dt$. One readily checks that the following decomposition holds: $\eta_i^n = {(\eta')}_i^n +{(\eta'')}_i^n$, with
\[
        (\eta')_i^n = \big(\overline{\Delta}_i^n(\drift{}{3}-\drift{}{2})\big)^2-\frac{\diffE{2}{3}{\vartheta}^2}{\diffE{1}{2}{\vartheta}^2} \big(\overline{\Delta}_i^n(\drift{}{2}-\drift{}{1})\big)^2
\]
and
\[
        (\eta'')_i^n=2\diffE{2}{3}{\vartheta}\Big(\overline{\Delta_i^n}\big((\drift{}{3}-\drift{}{2})-\frac{\diffE{2}{3}{\vartheta}}{\diffE{1}{2}{\vartheta}}(\drift{}{2}-\drift{}{1})\big)\Big)\int_{(i-1)\Delta_n}^{i\Delta_n} e^{\vartheta t}\sigmaS{t} d\Ws{t}.
\]
We will need the following lemma, proof of which is relatively straightforward yet technical and given in Section~\ref{subsection:preuveLemmeTechnique1}.
\begin{lem} \label{technical lemma 1}
        Let $(Y_t)_{t\geq 0}$ and $(Z_t)_{t \geq 0}$ be two c\`adl\`ag and progressively measurable processes. Assume that for some $s>1/2$, we have $\sup_{t \in [0,T]} t^{-s}\omega(Y)_t <\infty$. Then
        \[
                \Delta_n^{-1}\sum_{i = 1}^n \big(\overline{\Delta}_i^n Y\big)^2 \rightarrow \int_0^T Y_t^2dt
        \]
        and
        \[
                \Delta_n^{-1}\sum_{i = 1}^n\overline{\Delta}_i^n(Y)\int_{(i-1)\Delta_n}^{i\Delta_n}Z_tdB_t \rightarrow \int_0^T Y_tZ_tdB_t
        \]
        in probability.
\end{lem}
We successively have
\[
        \Delta_n^{-1}\sum_{i=1}^n(\eta')_i^n \rightarrow \int_0^T\mu_\vartheta(\drift{t}{})dt
\]
with $\mu_\vartheta(\drift{t}{})=(\drift{t}{3}-\drift{t}{2})^2-\frac{\diffE{2}{3}{\vartheta}^2}{\diffE{1}{2}{\vartheta}^2}(\drift{t}{2}-\drift{t}{1})^2$
and
\[
        \Delta_n^{-1}\sum_{i=1}^n(\eta'')_i^n \rightarrow 2\int_0^T\lambda_\vartheta(\drift{t}{})e^{\vartheta t}\sigmaS{t} d\Ws{t}
\]
in probability, by Lemma~\ref{technical lemma 1} applied to $Y_t=(\drift{t}{3}-\drift{t}{2})-\frac{\diffE{2}{3}{\vartheta}}{\diffE{1}{2}{\vartheta}}(\drift{t}{2}-\drift{t}{1})$ and $Z_t=e^{\vartheta t}\sigmaS{t}$, and Assumption \ref{asm:regulariteFonctions}, where $\lambda_\vartheta(\drift{t}{}) = \diffE{2}{3}{\vartheta}Y_t$.
This, together with \eqref{quad bis}, implies the convergence
\[
        \Delta_n^{-1}\big(\Psi_{T_1,T_2,T_3}^n -\psi_{T_1,T_2,T_3}(\vartheta)\big) \rightarrow \frac{\int_0^T\mu_\vartheta(\drift{t}{})dt+2\int_0^T\lambda_\vartheta(\drift{t}{})e^{\vartheta t}\sigmaS{t} d\Ws{t}}{\diffE{1}{2}{\vartheta}^2\int_0^Te^{2\vartheta t}\sigmaS{t}^2 dt}
\]
in probability.\par \bigskip

\noindent {\bf Step 3.} Finally, we have
\[
        \Delta_n^{-1}\big(\widehat{\vartheta}_{3,n}-\vartheta\big) = \Delta_n^{-1}\big(\Psi_{T_1,T_2,T_3}^n -\psi_{T_1,T_2,T_3}(\vartheta)\big)\partial_\vartheta\psi_{T_1,T_2,T_3}^{-1}(Z_n),
\]
for some $Z_n$ that converges to $\psi_{T_1,T_2,T_3}(\vartheta)$ by Step~1. Hence
\[
        \Delta_n^{-1}\big(\widehat{\vartheta}_{3,n}-\vartheta\big) \rightarrow \frac{\int_0^T\mu_\vartheta(\drift{t}{})dt+2\int_0^T\lambda_\vartheta(\drift{t}{})e^{\vartheta t}\sigmaS{t} d\Ws{t}}{\partial_\vartheta\psi_{T_1,T_2,T_3}(\vartheta)\diffE{1}{2}{\vartheta}^2\int_0^Te^{2\vartheta t}\sigmaS{t}^2 dt}
\]
and we conclude by noting that $\partial_\vartheta\psi_{T_1,T_2,T_3}(\vartheta)=\frac{2 D_3}{\diffE{1}{2}{\vartheta}^3}$.

\subsection{Proof of Proposition~\ref{thm:estimationNPunePeriode}} \label{subsection:preuveEstimationNPsimple}

We may assume that $b^j=0$, the case $b^j \neq 0$ being obtained exactly in the same line as in Step 4 of the proof of Theorem \ref{thm:loiAsymptotiqueHatAMoinsA} (1). For ease of notation, we write $\widehat{\vartheta}_{2,n}$ for $\max(\widehat{\vartheta}_{2,n},\varpi_n)$ and set $t_i=i\Delta_n$ for $i=1,\ldots, n$. We also define $K(t) = \indiq_{(0,1]}(t)$ and $K_h(t)=h^{-1}K(th^{-1})$ for $h>0$.
We have
\[
	\widehat{\sigmaS{}}_{n,t}^2-\sigmaS{t}^2 
	=\frac{\sum_{i=1}^n K_{h_n}\big(t-t_{i-1}\big)\big((\Delta_i^n X^1)^2-(\Delta_i^n X^2)^2\big)}{e^{-2\widehat{\vartheta}_{2,n}(T_1-t)}-e^{-2\widehat{\vartheta}_{2,n}(T_2-t)}}-\sigmaS{t}^2 = I+II,
\]
with
\begin{align*}
        I=&\Big(\frac{1}{e^{-2\widehat{\vartheta}_{2,n}(T_1-t)}-e^{-2\widehat{\vartheta}_{2,n}(T_2-t)}}-\frac{1}{e^{-2\vartheta(T_1-t)}-e^{-2\vartheta(T_2-t)}}\Big) \\
        &\times\sum_{i=1}^n K_{h_n}\big(t-t_{i-1}\big)\big((\Delta_i^n X^1)^2-(\Delta_i^n X^2)^2\big)
\end{align*}
and
\[
        II = \frac{\sum_{i=1}^n K_{h_n}\big(t-t_{i-1}\big)\big((\Delta_i^n X^1)^2-(\Delta_i^n X^2)^2\big)}{e^{-2\vartheta(T_1-t)}-e^{-2\vartheta(T_2-t)}}-\sigmaS{t}^2.
\]
\noindent {\bf Step 1.} Since $\E[\big(\Delta_i^n X^j\big)^2]$ is of order $\Delta_n$ by Burkholder-Davis-Gundy inequality, we have that $\E[\big|(\Delta_i^n X^1)^2-(\Delta_i^n X^2)^2\big|]$ is of order $\Delta_n$ as well and therefore  
\begin{align*}
         \E\Big[\Big|\sum_{i=1}^n K_{h_n}(t-t_{i-1})\big((\Delta_i^n X^1)^2-(\Delta_i^n X^2)^2\big)\Big|\Big] 
        \lesssim  & \sum_{i = 1}^nK_{h_n}(t-t_{i-1})\Delta_n \lesssim 1
\end{align*}
since $K_{h_n}(t-t_{i-1})$ is of order $h_n^{-1}$ for a number of terms that are at most of order $\Delta_n^{-1}h_n$. Therefore $ \sum_{i=1}^n K_{h_n}\big(t-t_{i-1}\big)\big((\Delta_i^n X^1)^2-(\Delta_i^n X^2)^2\big)$ is tight, and we conclude that $I$ is of order $\Delta_n^{1/2}$ in probability by applying Theorem \ref{thm:loiAsymptotiqueHatAMoinsA}~(1).\par \bigskip

\noindent {\bf Step 2.} The term $II$ further splits into $II = (e^{-2\vartheta(T_1-t)}-e^{-2\vartheta(T_2-t)})^{-1}\big(B_n(t)+V_n(t)\big)$,
having
\[
        V_n(t)=\sum_{i=1}^n K_{h_n}(t-t_{i-1})\big((\Delta_i^n X^1)^2-(\Delta_i^n X^2)^2-\E\big[(\Delta_i^n X^1)^2-(\Delta_i^n X^2)^2\big|\cF_{i-1}\big]\big)
\]
and
\[
        B_n(t)=\sum_{i=1}^n \E\big[K_{h_n}(t-t_{i-1})\big((\Delta_i^n X^1)^2-(\Delta_i^n X^2)^2\big)\big|\cF_{i-1}\big]-\big(e^{-2\vartheta(T_1-t)}-e^{-2\vartheta(T_2-t)}\big)\sigmaS{t}^2.
\]
Hereafter, we abbreviate $\cF_{i\Delta_n}$ by $\cF_{i}$.\\

\noindent {\bf Step 3.}
We first prove an upper bound for $\E[V_n(t)^2]$. We have
\begin{align*}
        & \sup_{t \in[h_n,T]} \E\Big[\Big(\sum_{i=1}^n K_{h_n}(t-t_{i-1})\big((\Delta_i^n X^1)^2-(\Delta_i^n X^2)^2-\E\big[(\Delta_i^n X^1)^2-(\Delta_i^n X^2)^2\big|\cF_{i-1}\big]\big)\Big)^2\Big] \\
        =& \sup_{t \in[h_n,T]}h_n^{-2} \sum_{i=1}^n K\Big(\frac{t-t_{i-1}}{h_n}\Big)^2 \E\Big[\big((\Delta_i^n X^1)^2-(\Delta_i^n X^2)^2-\E\big[(\Delta_i^n X^1)^2-(\Delta_i^n X^2)^2\big|\cF_{i-1}\big]\big)^2\Big] 
\end{align*}
because cross-terms in the development are zero due to conditioning. By compactness of the support of $K$, there are at most of order $\Delta_n^{-1}h_n$ nonvanishing terms in the sum and the estimate is uniform in $t \in[h_n,T]$. Finally, since
\[
        \E\Big[\big((\Delta_i^n X^1)^2-(\Delta_i^n X^2)^2-\E\big((\Delta_i^n X^1)^2-(\Delta_i^n X^2)^2\big|\cF_{i-1}\big)\big)^2\Big] \lesssim \Delta_n^2,
\]
we obtain $\sup_{t \in[h_n,T]}\E\big[V_n(t)^2\big]  \lesssim \Delta_n h_n^{-1}$.\\

\noindent {\bf Step 4.}
In order to bound the bias
we use the decomposition
\[
        B_n(t) = \big(e^{-2\vartheta(T_1-t)}-e^{-2\vartheta(T_2-t)}\big)(III+IV)),
\]
where
\[
        III=\int_0^Th_n^{-1}K\Big(\frac{t-u}{h_n}\Big)e^{-2\vartheta(t-u)}\sigmaS{u}^2du-\sigmaS{t}^2
\]
and
\begin{align*}
	IV=&\frac{\sum_{i=1}^n \E\Big[h_n^{-1}K\Big(\frac{t-t_{i-1}}{h_n}\Big)\big((\Delta_i^n X^1)^2-(\Delta_i^n X^2)^2\big)\Big|\cF_{i-1}\Big]}{e^{-2\vartheta(T_1-t)}-e^{-2\vartheta(T_2-t)}}- \int_0^Th_n^{-1}K\Big(\frac{t-u}{h_n}\Big)e^{-2\vartheta(t-u)}\sigmaS{u}^2du \text{.}
\end{align*}
For every $t \in[h_n,T]$ we have $\int_{\frac{t-T}{h_n}}^{\frac{t}{h_n}}K(x)dx=1$ hence
\begin{align*}
         \E\big[III^2\big]
        &= \E\Big[\Big(\int_{\frac{t-T}{h_n}}^{\frac{t}{h_n}}K(x)e^{-2\vartheta h_nx}\sigmaS{t-h_nx}^2dx-\sigmaS{t}^2\Big)^2 \Big] \\
       &=\E\Big[\Big(\int_{\textrm{supp}(K)}K(x)\big(e^{-2\vartheta h_nx}\sigmaS{t-h_nx}^2-\sigmaS{t}^2\big)dx\Big)^2 \Big] \\
        &\leq \int_{\textrm{supp}(K)}K(x)^2\E\big[\big(e^{-2\vartheta h_nx}\sigmaS{t-h_nx}^2-\sigmaS{t}^2\big)^2\big]dx
\end{align*}
by Jensen inequality since $\textrm{supp}(K)\subset\big[\frac{t-T}{h},\frac{t}{h}\big]$ and since $\int_0^T K(s)ds = 1$. By convexity,
\[
        (e^{-2\vartheta h_nx}\sigmaS{t-h_nx}^2-\sigmaS{t}^2)^2\leq 2(e^{-2\vartheta h_nx}\sigmaS{t-h_nx}^2-\sigmaS{t-h_nx}^2)^2+2(\sigmaS{t-h_nx}^2-\sigmaS{t}^2)^2
       \]
       follows. Bounding further the remainder in the expansion of $x\leadsto e^{-2\vartheta h_nx}$ at the point $0$, we obtain $|e^{-2\vartheta h_nx}-1|\leq M |2\vartheta h_nx|$ for some $M>0$. By localization, we find some $M_\sigma>0$ such that $\sigmaS{t}<M_\sigma$. It follows that 
\begin{align*}
        \E\big[III^2\big]&\leq \int_{\textrm{supp}(K)}K^2(x)\big(2M_\sigma^4(2\vartheta h_nx M)^2 +2\E\big(\big(\sigmaS{t-h_nx}^2-\sigmaS{t}^2\big)^2\big)\big)dx \\
        &\leq \int_{\textrm{supp}(K)}K^2(x)\big(2M_\sigma^4 (2\vartheta h_nx M)^2+2c|h_nx|^{2\alpha}\big)dx,
\end{align*}
using Assumption~\ref{asm:sigmasHolder}. This estimate is uniform in $t \in[h_n,T]$, therefore $\sup_{t \in[h_n,T]} \E\big[III^2\big] \lesssim h_n^{2\alpha}$.
Next, we write $IV =\sum_{i=1}^n h_n^{-1} \delta_i(t)$,
with
\[
        \delta_i(t)=\E\Big(K\Big(\frac{t-t_{i-1}}{h_n}\Big)\int_{t_{i-1}}^{t_i}e^{-2\vartheta(t-u)}\sigmaS{u}^2du\Big|\cF_{i-1}\Big)-\int_{t_{i-1}}^{t_i}K\Big(\frac{t-u}{h_n}\Big)e^{-2\vartheta (t-u)}\sigmaS{u}^2du \text{,}
\]
since the term
\begin{align*}
      & \sum_{i=1}^n h_n^{-1}K\Big(\frac{t-t_{i-1}}{h_n}\Big) \E\Big[\frac{\big((\Delta_i^n X^1)^2-(\Delta_i^n X^2)^2\big)}{e^{-2\vartheta(T_1-t)}-e^{-2\vartheta(T_2-t)}}-\int_{t_{i-1}}^{t_i}e^{-2\vartheta(t-u)}\sigmaS{u}^2du\Big|\cF_{i-1}\Big]\\
        &=\sum_{i=1}^n h_n^{-1}K\Big(\frac{t-t_{i-1}}{h_n}\Big) e^{-2\vartheta t} \E\Big[\Big(\int_{t_{i-1}}^{t_i} e^{\vartheta u}\sigmaS{u}d\Ws{u}\Big)^2-\int_{t_{i-1}}^{t_i}e^{2\vartheta u}\sigmaS{u}^2du\Big|\cF_{i-1}\Big] = 0
\end{align*}
vanishes. Therefore, writing
$\E[IV^2] = \sum_{i=1}^n \E\big(h_n^{-2}(\delta_i(t))^2\big) +2\sum_{1\leq i<j\leq n} \E\big[h_n^{-2} \delta_i(t) \delta_j(t) \big]$, we successively obtain
\[
        \E[h_n^{-2}(\delta_i^\cS(t))^2]\leq \frac{\Delta_n}{h_n^2}\int_{t_{i-1}}^{t_i} e^{-4\vartheta (t-u)} \E\big[\big(K\big(\frac{t-t_{i-1}}{h_n}\big)\E[\sigmaS{u}^2|\cF_{i-1}]-K\big(\frac{t-u}{h_n}\big)\sigmaS{u}^2\big)^2\big]
\]
by Jensen inequality, so that $\sum_{i=1}^n \E\big[h_n^{-2}(\delta_i^\cS(t))^2\big]\lesssim \Delta_n h_n^{-1}$ uniformly in $t$, since there are at most of order $\Delta_n h_n$ nonvanishing terms in the sum. Finally, conditioning on $\cF_{j-1}$,
\[
        \E\big[h_n^{-2} \delta_i(t) \delta_j(t) \big]=h_n^{-2}\E\Big[\delta_i(t)\E\big(\int_{t_{j-1}}^{t_j}\big(K\big(\frac{t-t_{j-1}}{h_n}\big)-K\big(\frac{t-u}{h_n}\big)\big)e^{-2\vartheta(t-u)}\sigmaS{u}^2 du\big|\cF_{j-1}\big]\Big],
\]
and the difference $K\Big(\frac{t-t_{j-1}}{h_n}\Big)-K\Big(\frac{t-u}{h_n}\Big)$ is non-zero only if $t\in(t_{j-1},u]$ or $t\in(t_{j-1}+h_n,u+h_n]$, which can be the case for $j$ in some set $\cJ_t$, which contains at most three indexes. Therefore,
\begin{align*}
        \Big|\sum_{1\leq i<j\leq n} \E\big[h_n^{-2} \delta_i(t) \delta_j(t) \big]\Big|&=\Big|\sum_{i=1}^{n-1} \sum_{j\in\cJ_t} \E\big[h_n^{-2} \delta_i(t) \delta_j(t) \big]\Big|
        \leq 3h_n^{-2}\sum_{i=1}^{n-1} \E\big[\big|\delta_i(t)\big|\big]M_\sigma^2e^{2\vartheta T} \Delta_n,
\end{align*}
which is of order $\Delta_n h_n^{-1}$. We infer $\sup_{t \in[h_n,T]} \E\big[IV^2\big] \lesssim \Delta_n h_n^{-1}$.\\

\noindent {\bf Step 5.} From the estimates established in Steps 3. and 4. we derive
$$\sup_{t \in[h_n,T]}\E\big[(V_n(t))^2\big] \lesssim \Delta_n h_n^{-1},\;\sup_{t \in[h_n,T]} \E\big[B_n(t)^2\big] \lesssim h_n^{2\alpha}+\Delta_n h_n^{-1}.$$
The choice $h_n = \Delta_n^{1/(2\alpha+1)}$ implies that the two error terms $h_n^{2\alpha}$ and $\Delta_nh_n^{-1}$ are of the same order, namely $\Delta_n^{2\alpha/(2\alpha+1)}$, 
which ends the proof concerning $\sigmaS{}$.\\

\noindent {\bf Step 6.} The proof is the same for $\widehat{\sigmaL{}}^2$. We split $\widehat{\sigmaL{}}_{n,t}^2-\sigmaL{t}^2$ as follows
\begin{align*}
	\widehat{\sigmaL{}}_{n,t}^2-\sigmaL{t}^2= & \Big(\frac{e^{-2\widehat{\vartheta}_{2,n}T_2}}{e^{-2\widehat{\vartheta}_{2,n}T_1}-e^{-2\widehat{\vartheta}_{2,n}T_2}}-\frac{e^{-2\vartheta T_2}}{e^{-2\vartheta T_1}-e^{-2\vartheta T_2}}\Big)\sum_{i=1}^n h_n^{-1}K\Big(\frac{t-t_{i-1}}{h_n}\Big)(\Delta_i^n X^1)^2 \\
	&+\Big(\frac{e^{-2\widehat{\vartheta}_{2,n}T_1}}{e^{-2\widehat{\vartheta}_{2,n}T_1}-e^{-2\widehat{\vartheta}_{2,n}T_2}}-\frac{e^{-2\vartheta T_1}}{e^{-2\vartheta T_1}-e^{-2\vartheta T_2}}\Big)\sum_{i=1}^n h_n^{-1}K\Big(\frac{t-t_{i-1}}{h_n}\Big)(\Delta_i^n X^2)^2 \\
	&+\frac{\sum_{i=1}^n h_n^{-1}K\Big(\frac{t-t_{i-1}}{h_n}\Big)(e^{-2\vartheta (T_1-t)}(\Delta_i^n X^2)^2-e^{-2\vartheta (T_2-t)}(\Delta_i^n X^1)^2)}{e^{-2\vartheta (T_1-t)}-e^{-2\vartheta (T_2-t)}}-\sigmaL{t}^2 
\end{align*}
and proceed analogously. 
The proof of Theorem \ref{thm:estimationNPunePeriode} is complete.

\subsection{Proof of Theorem~\ref{thm:borneInferieureVarianceEstimationTheta}} \label{subsection:preuveTheoremeBorneInf}

\subsubsection*{Preliminaries on efficient semiparametric estimation}
We refer to Sections~25.3--25.4 of \cite{bib:vanDerVaartAsymptoticStatistics} for a comprehensive presentation of efficient semiparametric estimation, that we need to adapt to our framework.
Assuming $(\sigmaS,\sigmaL{})$ to be deterministic, the data $(\Delta_i^n X^1, \Delta_i^n X^2), i=1,\ldots,n$ extracted from \eqref{diff gaussienne} generate a product experiment $\mathcal E^n = \otimes_{i = 1}^n \mathcal P^{i,n}$, where
\[
        \cP^{i,n}=\big\{f_{\vartheta,\sigmaS{},\sigmaL{}}^{i,n}, (\vartheta, \sigmaS{},\sigmaL{})\in [0,\infty) \times \Sigma(c,\tilde{c}\big\},
\]
where $ f_{\vartheta,\sigmaS{},\sigmaL{}}^{i,n}$ is the density on $\R^2$ of the Gaussian vector $(\Delta_i^n X^1, \Delta_i^n X^2)$, see Section \ref{subsection:estimationEfficaceTheta}.
Let $\varepsilon >0$ considered as a small perturbation parameter.  For every $u \in [0,\varepsilon]$, let us be given moreover two regular functions $t \leadsto \etaS{}^u_t$ and $t \leadsto \etaL{}^u_t$ such that
\begin{equation} \label{def perturbation}
\tilde c < \inf_{t \in [0,T],u \in [0,\varepsilon]} \min(\etaS{}^u_t, \etaL{}^u_t) \leq \sup_{t \in [0,T],u \in [0,\varepsilon]}\max(\etaS{}^u_t, \etaL{}^u_t) < c.
\end{equation}
We call  $\big(t \leadsto (\etaS{}^u_t, \etaL{}^u_t)\big)_{0 \leq u \leq \varepsilon}$ a perturbation path if we have \eqref{def perturbation} together with the property
$$(\etaS{t}^0, \etaL{t}^0)= (\sigmaS{t},\sigmaL{t})\;\;\text{for every}\;\;t \in [0,T].$$
Let $\vartheta$ and $\iota \in \R$ be such that $\vartheta+\iota \varepsilon>0$.  Remember that $ f_{\vartheta,\sigmaS{},\sigmaL{}}^{i,n}$ is the density on $\R^2$ of the Gaussian vector $(\Delta_i^n X^1, \Delta_i^n X^2)$ extracted from \eqref{diff gaussienne}. We obtain a parametric submodel of $\mathcal P^{i,n}$ around $\big(t\leadsto (\sigmaS{t}^2,\sigmaL{t}^2)\big)$ by setting 
$$\mathcal P_0^{i,n} =\big\{ f_{\vartheta+\iota u,\etaS{}^u,\etaL{}^u}^{i,n}, 0\leq u \leq \varepsilon \big\},$$
noting that $\cP_0^{i,n}$ passes through the true distribution at $u=0$.
We consider only submodels that are differentiable in quadratic mean at $u=0$, with score function $g_{i, \iota, \etaS{},\etaL{}}\in L^2(\P_{\vartheta,\sigmaS{},\sigmaL{}})$. If we let $\cP_0^{i,n}$ range over all admissible submodels as $(\eta, \overline{\eta})$ varies among perturbation paths, we obtain a collection of score functions that define in turn the tangent set $\dot{\cP}_{\vartheta, \sigmaS{},\sigmaL{}}^{i,n}$ of the model $\cP^{i,n}$ at the true distribution.
Any score function $g_{\iota, \etaS{},\etaL{}}^{i,n} \in \dot{\cP}_{\vartheta, \sigmaS{},\sigmaL{}}^{i,n}$ admits the representation
\begin{equation} \label{eq: rep score}
g_{\iota, \etaS{},\etaL{}}^{i,n}=\iota \ell_{\vartheta, \sigmaS{},\sigmaL{}}^{i,n}+g_{\etaS{}, \etaL{}}^{i,n},
\end{equation}
where $\ell_{\vartheta, \sigmaS{},\sigmaL{}}^{i,n}$ is the score function of the original model defined in \eqref{eq: def score} when $\sigmaS{}$ and $\sigmaL{}$ are known, and $g_{\etaS{},\etaL{}}^{i,n}$ is the score function obtained from a parametric submodel at $\iota=0$, to be interpreted as the score relative to the nuisance parameter, while $\ell_{\vartheta, \sigmaS{},\sigmaL{}}^{i,n}$ corresponds to the score relative to the parameter of interest $\vartheta$.
\subsubsection*{Completion of proof of Theorem~\ref{thm:borneInferieureVarianceEstimationTheta}} 

From \eqref{eq: def score}, and the explicit representation 
  $$  f_{\vartheta,\sigmaS{},\sigmaL{}}^{i,n}(x,y)  =  \frac{\exp\Big(-\frac{\big(\int_{(i-1)\Delta_n}^{i\Delta_n}\sigmaL{t}^2dt\big)(x-y)^2+\big(\int_{(i-1)\Delta_n}^{i\Delta_n}e^{-2\vartheta(T_1-t)}\sigmaS{t}^2dt\big)\big(y-e^{-\vartheta(T_2-T_1)}x\big)^2}{2\big(\int_{(i-1)\Delta_n}^{i\Delta_n}e^{-2\vartheta(T_1-t)}\sigmaS{t}^2 dt\big) \big(\int_{(i-1)\Delta_n}^{i\Delta_n}\sigmaL{t}^2dt\big)\big(1-e^{-\vartheta(T_2-T_1)}\big)^2}\Big)}{2\pi \big(\int_{(i-1)\Delta_n}^{i\Delta_n}e^{-2\vartheta(T_1-t)}\sigmaS{t}^2dt\big)^{1/2} \big(\int_{(i-1)\Delta_n}^{i\Delta_n}\sigmaL{t}^2dt\big)^{1/2}\big(1-e^{-\vartheta(T_2-T_1)}\big)},
  $$
we derive
\begin{align*}
     \partial_{\vartheta} \log f_{\vartheta,\sigmaS{},\sigmaL{}}^{i,n}(\Delta_i^n X^1,\Delta_i^n X^2) & = \frac{\int_{(i-1)\Delta_n}^{i\Delta_n}(T_1-t)e^{-2\vartheta(T_1-t)}\sigmaS{t}^2dt}{\int_{(i-1)\Delta_n}^{i\Delta_n}e^{-2\vartheta(T_1-t)}\sigmaS{t}^2dt}-\frac{e^{-\vartheta(T_2-T_1)}(T_2-T_1)}{1-e^{-\vartheta(T_2-T_1)}} \\
        &+ (\Delta_i^n X^2-\Delta_i^n X^1)^2 \frac{\int_{(i-1)\Delta_n}^{i\Delta_n}(e^{-\vartheta(T_2-T_1)}(T_2-t)-(T_1-t))e^{-2\vartheta(T_1-t)}\sigmaS{t}^2dt}{\big(\int_{(i-1)\Delta_n}^{i\Delta_n}e^{-2\vartheta(T_1-t)}\sigmaS{t}^2dt)^2(1-e^{-\vartheta(T_2-T_1)})^3} \\
        &+ (\Delta_i^n X^2-e^{-\vartheta(T_2-T_1)}\Delta_i^n X^1)^2 \frac{(T_2-T_1)e^{-\vartheta(T_2-T_1)}}{\int_{(i-1)\Delta_n}^{i\Delta_n}\sigmaL{t}^2dt (1-e^{-\vartheta(T_2-T_1)})^3} \\
        &- \Delta_i^n X^1(\Delta_i^n X^2-e^{-\vartheta(T_2-T_1)}\Delta_i^n X^1) \frac{e^{-\vartheta(T_2-T_1)}(T_2-T_1)}{\int_{(i-1)\Delta_n}^{i\Delta_n}\sigmaL{t}^2dt (1-e^{-\vartheta(T_2-T_1)})^2}.
\end{align*}
We pick a path $(\etaS{t}^u, \etaL{t}^u)_{0 \leq u \leq \varepsilon}$ of the form
$\etaS{t}^u=(1+uk(t))\sigmaS{t}$ and $\etaL{t}^u=(1+u\overline{k}(t))\sigmaL{t}$ so that $(\etaS{}^u,\etaL{}^u)\in\Sigma(c,\tilde{c})$.
The submodel is differentiable in quadratic mean at $u=0$, with score function having representation
$$g_{\iota, k, \overline{k}}^{i,n} = \iota \ell_{\vartheta,\sigmaS{},\sigmaL{}}^i + g_{k, \overline{k}}^{i,n}.$$
according to \eqref{eq: rep score} and parametrised by $(k,\overline{k})$.
Formally, 
$$\iota \ell_{\vartheta,\sigmaS{},\sigmaL{}}^{i,n} = \big(\tfrac{d}{du}\log f_{\vartheta+\iota  u,\sigmaS{},\sigmaL{}}^{i,n}\big)_{u = 0}\;\;\text{and}\;\;g_{k, \overline{k}}^{i,n}=\big(\tfrac{d}{du}\log f_{\vartheta,\etaS{}^u,\etaL{}^u}^{i,n}\big)_{u=0},$$ 
so that
\begin{align*}
        g_{k, \overline{k}}^{i,n}=& -\frac{\int_{(i-1)\Delta_n}^{i\Delta_n}e^{-2\vartheta(T_1-t)}\sigmaS{t}^2k(t)dt}{\int_{(i-1)\Delta_n}^{i\Delta_n}e^{-2\vartheta(T_1-t)}\sigmaS{t}^2dt} - \frac{\int_{(i-1)\Delta_n}^{i\Delta_n}\sigmaL{t}^2\overline{k}(t)dt}{\int_{(i-1)\Delta_n}^{i\Delta_n}\sigmaL{t}^2dt} \\
        &+ (\Delta_i^n X^2-\Delta_i^n X^1)^2 \frac{\int_{(i-1)\Delta_n}^{i\Delta_n}e^{-2\vartheta(T_1-t)}\sigmaS{t}^2k(t)dt}{\big(\int_{(i-1)\Delta_n}^{i\Delta_n}e^{-2\vartheta(T_1-t)}\sigmaS{t}^2dt\big)^2(1-e^{-\vartheta(T_2-T_1)})^2} \\
        &+ (\Delta_i^n X^2-e^{-\vartheta(T_2-T_1)}\Delta_i^n X^1)^2 \frac{\int_{(i-1)\Delta_n}^{i\Delta_n}\sigmaL{t}^2\overline{k}(t)dt}{\big(\int_{(i-1)\Delta_n}^{i\Delta_n}\sigmaL{t}^2dt\big)^2(1-e^{-\vartheta(T_2-T_1)})^2} \text{.}
\end{align*}
Introduce the orthogonal projection $\Pi$ onto (the closure of) of $\mathrm{Span}\{g_{k, \overline{k}}^{i,n}, \text{for all admissible}\;(k,\overline{k})\}$. 
Then $\widetilde{\ell}_{\vartheta, \sigmaS{},\sigmaL{}}^{i,n}=\ell_{\vartheta, \sigmaS{},\sigmaL{}}^{i,n}-\Pi \ell_{\vartheta, \sigmaS{},\sigmaL{}}^{i,n}$ is the efficient score for $\vartheta$ and $\widetilde{I}_{\vartheta, \sigmaS{},\sigmaL{}}^{i,n}=\E_{\vartheta,\sigmaS{},\sigmaL{}}\big[(\widetilde{\ell}_{\vartheta, \sigmaS{},\sigmaL{}}^{i,n})^2\big]$ is the best achievable information bound, see Sections~25.3--25.4 of \cite{bib:vanDerVaartAsymptoticStatistics} for details. By orthogonality, 
\[
      \E_{\vartheta,\sigmaS{},\sigmaL{}}\big[(\ell_{\vartheta, \sigmaS{},\sigmaL{}}^{i,n} -\Pi \ell_{\vartheta, \sigmaS{},\sigmaL{}}^{i,n})g_{k, \overline{k}}^{i,n}\big]=0,\;\;\text{for all admissible}\;\;(k,\overline{k}),
\]
and anticipating further the representation $\Pi \ell_{\vartheta, \sigmaS{},\sigmaL{}}^{i,n} = g_{k^\star, \overline{k}^\star}^{i,n}$ for some admissible $(k^\star,\overline{k}^\star)$, it suffices to solve
\[
        0=\E_{\vartheta,\sigmaS{},\sigmaL{}}\big[\big(l_{\vartheta, \sigmaS{},\sigmaL{}}^{i,n} -g_{k^\star, \overline{k}^\star}^{i,n}\big)g_{k, \overline{k}}^{i,n}\big]\;\;\text{ for all admissible}\;\;(k,\overline{k}).
\]
Elementary computations yield $k^\star(t)=\frac{(T_2-t)e^{-\vartheta (T_2-T_1)}t(T_1-t)}{1-e^{-\vartheta(T_2-T_1)}}$ and $\overline{k}^\star(t)=0$. We conclude
\[
   \widetilde{\ell}_{\vartheta, \sigmaS{},\sigmaL{}}^{i,n}  = \widetilde{\ell}_{\vartheta,\sigmaL{}}^{i,n}    
    =\frac{(\Delta_i^n X^2-\Delta_i^n X^1)(\Delta_i^n X^2-e^{-\vartheta(T_2-T_1)}\Delta_i^n X^1)e^{-\vartheta(T_2-T_1)}(T_2-T_1)}{(1-e^{-\vartheta(T_2-T_1)})^3\int_{(i-1)\Delta_n}^{i\Delta_n}\sigmaL{t}^2dt}
\]
and is independent of $\sigma$. Furthermore, the best achievable information bound becomes
\[
        \widetilde{I}_{\vartheta, \sigmaS{},\sigmaL{}}^{i,n}=\frac{(T_2-T_1)^2}{(e^{\vartheta(T_2-T_1)}-1)^2}\frac{\int_{(i-1)\Delta_n}^{i\Delta_n}e^{-2\vartheta(T_1-t)}\sigmaS{t}^2dt}{\int_{(i-1)\Delta_n}^{i\Delta_n}\sigmaL{t}^2dt} \text{.}
\]
By independence of the increments $(\Delta_i^n X^1, \Delta_i^n X^2)$, it remains to piece together the results for each $\mathcal P^{i,n}$ using the product structure of $\mathcal E^n$.
We find the asymptotically equivalent bound 
\[
    \sum_{i = 1}^n \widetilde{I}_{\vartheta, \sigmaS{},\sigmaL{}}^{i,n} \sim n\frac{(T_2-T_1)^2}{T(e^{\vartheta(T_2-T_1)}-1)^2}\int_0^T\frac{e^{-2\vartheta(T_1-t)}\sigmaS{t}^2}{\sigmaL{t}^2}dt
\]
as $n \rightarrow \infty$. Taking the inverse and dividing by $\Delta_n$ we obtain the desired bound. The proof of Theorem \ref{thm:borneInferieureVarianceEstimationTheta} is complete.

\subsection{Proof of Theorem~\ref{thm:estimationEfficace}} \label{subsection:preuveEstimationEfficaceSimple}
The first assertion was obtained in Section~\ref{subsection:preuveTheoremeBorneInf} in the course of the proof of Theorem~\ref{thm:borneInferieureVarianceEstimationTheta}.
With no loss of generality, we may (and will) assume that $b^j=0$,  the case $b^j \neq 0$ being obtained exactly in the same line as in Step 4 of the proof of Theorem \ref{thm:loiAsymptotiqueHatAMoinsA} (1). We also assume with no loss of generality that $\sup_{t} \max(\sigma_t,\overline{\sigma}_t) \leq M$. 
We further abbreviate $i\Delta_n$ by $t_i$.\\

\noindent {\bf Step 1.} Let $\vartheta_n$ be a deterministic sequence such that $\sqrt{n}(\vartheta_n-\vartheta)$ is bounded. 
We first prove
\begin{equation} \label{eqn:premiereEtapePreuveEfficacite}
\Delta_n^{1/2}\sum_{i\in\cI_n} \big(\widetilde{\ell}_{\vartheta_n, \sigmaL{}}^{i,n}-\widetilde{\ell}_{\vartheta_n, \widehat{\sigmaL{}}_n}^{i,n}\big)\to 0
\end{equation}
in probability, as $n\to\infty$. We have
\begin{align*}
&\; \Delta_n^{1/2}\sum_{i\in\cI_n} \big(\widetilde{\ell}_{\vartheta_n, \sigmaL{}}^{i,n}-\widetilde{\ell}_{\vartheta_n, \widehat{\sigmaL{}}_n}^{i,n}\big) \\
= &\;   \Delta_n^{1/2}\sum_{i\in\cI_n} (\Delta_i^n X^2-\Delta_i^n X^1)(\Delta_i^n X^2-e^{-\vartheta_n (T_2-T_1)}\Delta_i^n X^1)\Big(\frac{1}{\int_{t_{i-1}}^{t_i}\sigmaL{t}^2dt}-\frac{1}{\Delta_n\widehat{\sigmaL{}}_{n,t_{i-1}}^2}\Big) \\
= &\; I +II,
\end{align*}
with
\begin{align*}
I &= \Delta_n^{-1/2}\sum_{i\in\cI_n} (\Delta_i^n X^2-\Delta_i^n X^1)(\Delta_i^n X^2-e^{-\vartheta_n (T_2-T_1)}\Delta_i^n X^1) \tfrac{\Delta_n\sigmaL{t_{i-1}}^2-\int_{t_{i-1}}^{t_i}\sigmaL{t}^2dt}{\sigmaL{t_{i-1}}^2\int_{t_{i-1}}^{t_i}\sigmaL{t}^2dt},\\
II& = \Delta_n^{-1/2}\sum_{i\in\cI_n} (\Delta_i^n X^2-\Delta_i^n X^1)(\Delta_i^n X^2-e^{-\vartheta_n (T_2-T_1)}\Delta_i^n X^1)\tfrac{\widehat{\sigmaL{}}_{n,t_{i-1}}^2-\sigmaL{t_{i-1}}^2}{\sigmaL{t_{i-1}}^2\widehat{\sigmaL{}}_{n,t_{i-1}}^2}.
\end{align*}
 
For the term $I$, by Cauchy-Schwarz inequality and the fact that $\P((\sigmaS{},\sigmaL{})\in\Sigma(c,\tilde{c}))=1$, we have
\begin{align*}
        &\E\big[\big|\Delta_n^{-1/2}(\Delta_i^n X^2-\Delta_i^n X^1)(\Delta_i^n X^2-e^{-\vartheta_n (T_2-T_1)}\Delta_i^n X^1) \tfrac{\Delta_n\sigmaL{t_{i-1}}^2-\int_{t_{i-1}}^{t_i}\sigmaL{t}^2dt}{\sigmaL{t_{i-1}}^2\int_{t_{i-1}}^{t_i}\sigmaL{t}^2dt}\big|\big]\\
        \leq &\,\Delta_n^{-3/2}\tilde{c}^{-4}\E\big[\big((\Delta_i^n X^2-\Delta_i^n X^1)(\Delta_i^n X^2-e^{-\vartheta_n (T_2-T_1)}\Delta_i^n X^1)\big)^2\big]^{1/2}\E\big[\big(\Delta_n\sigmaL{t_{i-1}}^2-\int_{t_{i-1}}^{t_i}\sigmaL{t}^2dt\big)^2\big]^{1/2}.
\end{align*}
Combining Cauchy-Schwarz and Burkholder-Davis-Gundy inequalities and the smoothness Assumption \ref{asm:sigmasHolder} we successively obtain
\begin{align*} 
&\E\big[\big((\Delta_i^n X^2-\Delta_i^n X^1)(\Delta_i^n X^2-e^{-\vartheta_n (T_2-T_1)}\Delta_i^n X^1)\big)^2\big]  \lesssim \Delta_n^2, \\
&\E\big[\big(\Delta_n\sigmaL{t_{i-1}}^2-\int_{t_{i-1}}^{t_i}\sigmaL{t}^2dt\big)^2\big]  \lesssim \Delta_n^{2(1+\alpha)}.
\end{align*}
We infer $\E[|I|] \lesssim \sum_{i\in\cI_n}\Delta_n^{-3/2}\Delta_n\Delta_n^{1+\alpha} \lesssim \Delta_n^{\alpha-1/2} \rightarrow 0$ since $\alpha >1/2$ by assumption. For the term $II$, since
the kernel $K$ used for the nonparametric estimation has support included in $[0,\infty)$, we have that $\widehat{\sigmaL{}}_{n,t_{i-1}}^2$ is $\cF_{i-1}$-measurable. Conditioning on $\cF_{i-1}$, we set
\begin{align*}
        \chi_i^n  & =    \E\big[\Delta_n^{-1/2} (\Delta_i^n X^2-\Delta_i^n X^1)(\Delta_i^n X^2-e^{-\vartheta_n (T_2-T_1)}\Delta_i^n X^1)\tfrac{\widehat{\sigmaL{}}_{n,t_{i-1}}^2-\sigmaL{t_{i-1}}^2}{\sigmaL{t_{i-1}}^2\widehat{\sigmaL{}}_{n,t_{i-1}}^2}\big|\cF_{i-1}\big]\\
   & =\Delta_n^{-1/2}(e^{-\vartheta(T_2-T_1)}-e^{-\vartheta_n(T_2-T_1)})(e^{-\vartheta(T_2-T_1)}-1)\xi_i^n,
\end{align*}
say, with $\xi_i^n=\E\big[\int_{t_{i-1}}^{t_i}e^{-2\vartheta(T_1-t)}\sigmaS{t}^2dt\big|\cF_{i-1}\big]\tfrac{\widehat{\sigmaL{}}_{n,t_{i-1}}^2-\sigmaL{t_{i-1}}^2}{\sigmaL{t_{i-1}}^2\widehat{\sigmaL{}}_{n,t_{i-1}}^2}$. It follows that
$$
       \sum_{i=1}^n  \E\big[|\xi_i^n|\big|\cF_{i-1}\big]\leq \Delta_n M^2\tilde{c}^{-4}\sum_{i=1}^n \sup_{i\in\cI_n}\big|\widehat{\sigmaL{}}_{n,t_{i-1}}^2-\sigmaL{t_{i-1}}^2\big| \rightarrow 0
$$
in probability by Theorem~\ref{thm:estimationNPunePeriode}. Since $\Delta_n^{-1/2}\big(e^{-\vartheta(T_2-T_1)}-e^{-\vartheta_n(T_2-T_1)}\big)$ is bounded, we use Lemma~3.4 in \cite{bib:statisticsAndHFDataInStatisticalMethodsForSDE} applied to variables $\xi_i^n$ to conclude
   $     \sum_{i=1}^{\lfloor t/ \Delta_n\rfloor}
        \E[\chi_i^n\,|\,\mathcal F_{i-1}] \rightarrow 0
        $
in probability, locally uniformly in $t$. Moreover,
\begin{align*}
        \E\big[\big(\chi_i^n\big)^2        \big|\cF_{i-1}\big] &
        \lesssim \Delta_n^{-1}\E\big[\big((\Delta_i^n X^2-\Delta_i^n X^1)(\Delta_i^n X^2-e^{-\vartheta_n (T_2-T_1)}\Delta_i^n X^1)\big)^2\big|\cF_{i-1}\big]\sup_{i\in\cI_n}|\widehat{\sigmaL{}}_{n,t_{i-1}}^2-\sigmaL{t_{i-1}}^2| 
\end{align*}
which is of order $\Delta_n^{-1}\Delta_n^2 \Delta_n^{\alpha/(2\alpha+1)}$ so that 
$\sum_{i\in\cI_n}
\E\big[\big(\chi_i^n\big)^2        \big|\cF_{i-1}\big] \lesssim \Delta_n^{\alpha/(2\alpha+1)} \rightarrow 0$
in probability. Applying Lemma~3.4 in \cite{bib:statisticsAndHFDataInStatisticalMethodsForSDE} to the sequence $\chi_i^n$ enables us to conclude that $II$ converges to $0$ in probability and \eqref{eqn:premiereEtapePreuveEfficacite} follows.\\

\noindent {\bf Step 2.}
Since $\vartheta\leadsto \widetilde{\ell}_{\vartheta,\sigmaL{}}^{i,n}$ is smooth (at least twice differentiable) and $\vartheta_n-\vartheta$ is of order $n^{-1/2}$ a second-order Taylor expansion of $\widetilde{\ell}_{\vartheta_n,\sigmaL{}}^{i,n}$ at $\vartheta$ implies that 
\begin{equation} \label{eq: Taylor}
\sqrt{n}\big(\widetilde \ell_{\vartheta_n,\sigmaL{}}^{i,n}-\widetilde \ell_{\vartheta,\sigmaL{}}^{i,n}\big) - \sqrt{n}(\vartheta_n-\vartheta) \partial_{\vartheta}\widetilde \ell_{\vartheta,\sigmaL{}}^{i,n} \rightarrow 0
\end{equation}
in probability. Since $|\cI_n|\sim n(1-h_n)\sim n$, it is not difficult to check that 
\begin{equation} \label{eq: equivalence Fisher}
\Delta_n\sum_{i\in\cI_n}  \partial_{\vartheta}\widetilde \ell_{\vartheta,\sigmaL{}}^{i,n} \sim - \widetilde I_{\vartheta,\sigma{},\sigmaL{}} = - \frac{(T_2-T_1)^2}{(e^{\vartheta T_2}-e^{\vartheta T_1})^2}\int_0^T \frac{e^{2\vartheta t}\sigmaS{t}^2}{\sigmaL{t}^2} dt
\end{equation}
in probability under $\mathbb P_{\vartheta,\sigma{},\sigmaL{}}$,
that is the total Fisher information associated to the efficient scores $\widetilde{\ell}_{\vartheta,\sigmaL{}}^{i,n}$. Summing each term in \eqref{eq: Taylor} for  $i \in \mathcal I_n$, using \eqref{eq: equivalence Fisher} and the fact that $\Delta_n^{-1}$ and $n$ are of the same order, we further infer
\begin{equation}   \label{eqn:dlAnA}
 \Delta_n^{1/2}\Big(\sum_{i\in\cI_n} \big(\widetilde{\ell}_{\vartheta_n, \widehat{\sigmaL{}}_n}^{i,n}-\widetilde{\ell}_{\vartheta,\sigmaL{}}^{i,n}\big)+\Delta_n^{-1} (\vartheta_n-\vartheta)\widetilde{I}_{\vartheta, \sigmaS{},\sigmaL{}}\Big) \rightarrow 0
      \end{equation}
in probability, using \eqref{eqn:premiereEtapePreuveEfficacite} in order to substitute $\sigmaL{}$ by $\widehat{\sigmaL{}}_n$ in the first term. Moreover, \eqref{eqn:dlAnA} remains true if we replace $\vartheta_n$ by the discretised version $\widehat{\vartheta}_{2,n}$ of $\hat \vartheta_n$, using moreover that $\sqrt{n}(\hat \vartheta_n-\vartheta)$ is bounded in probability thanks to Theorem \ref{thm:loiAsymptotiqueHatAMoinsA}. We refer to the proof of Theorem~5.48 in \cite{bib:vanDerVaartAsymptoticStatistics} for the details.\\
\noindent {\bf Step 3.} We establish
\begin{equation}
        \Delta_n\sum_{i\in\cI_n} \big(\widetilde{\ell}_{\widehat{\vartheta}_{2,n}, \widehat{\sigmaL{}}_n}^{i,n}\big)^2 \to \widetilde{I}_{\vartheta,\sigmaS{},\sigmaL{}} \label{eqn:convergenceVersInformationAvecEstimateursSigma}
\end{equation}
in probability under $\mathbb P_{\vartheta,\sigmaS{},\sigmaL{}}$ 
. The computations are similar to Step 1, combining Theorems \ref{thm:loiAsymptotiqueHatAMoinsA} and \ref{thm:estimationNPunePeriode}. We briefly give the mains steps. Observe that 
$$\Delta_n\sum_{i\in\cI_n} \big(\widetilde{\ell}_{\vartheta_n, \widehat{\sigmaL{}}_n}^{i,n}\big)^2-\Delta_n\sum_{i\in\cI_n} \big(\widetilde{\ell}_{\vartheta,\sigmaL{}}^{i,n}\big)^2 = (T_2-T_1)\sum_{i\in\cI_n}(I_i+II_i+III_i),$$
with
\begin{align*}
        I_i & =\Delta_n\tfrac{(\Delta_i^nX^2-\Delta_i^nX^1)^2}{\big(\Delta_n\widehat{\sigmaL{}}_{n,t_{i-1}}^2\big)^2}\big(\tfrac{(\Delta_i^n X^2-e^{-\vartheta_n (T_2-T_1)}\Delta_i^n X^1)^2 e^{-2\vartheta_n (T_2-T_1)}}{(1-e^{-\vartheta_n(T_2-T_1)})^6}- \tfrac{(\Delta_i^n X^2-e^{-\vartheta(T_2-T_1)}\Delta_i^n X^1)^2 e^{-2\vartheta(T_2-T_1)}}{(1-e^{-\vartheta(T_2-T_1)})^6}\big),\\
        II_i & = \Delta_n(\Delta_i^nX^2-\Delta_i^nX^1)^2 \tfrac{(\Delta_i^n X^2-e^{-\vartheta(T_2-T_1)}\Delta_i^n X^1)^2 e^{-2\vartheta(T_2-T_1)}}{(1-e^{-\vartheta(T_2-T_1)})^6}\big(\tfrac{1}{\big(\Delta_n\widehat{\sigmaL{}}_{n,t_{i-1}}^2\big)^2}-\tfrac{1}{\big(\Delta_n\sigmaL{t_{i-1}}^2\big)^2}\big),\\
        III_i & = \Delta_n(\Delta_i^nX^2-\Delta_i^nX^1)^2 \tfrac{(\Delta_i^n X^2-e^{-\vartheta(T_2-T_1)}\Delta_i^n X^1)^2 e^{-2\vartheta(T_2-T_1)}}{(1-e^{-\vartheta(T_2-T_1)})^6}\big(\tfrac{1}{\big(\Delta_n\sigmaL{t_{i-1}}^2\big)^2}-\frac{1}{\big(\int_{t_{i-1}}^{t_i}\sigmaL{t}^2dt\big)^2}\big).
\end{align*}
For the term $I_i$ we use a Taylor expansion of $\vartheta_n$ near $\vartheta$ in order to obtain
$\E[|\cT_{i,n}^1|]\lesssim |\vartheta_n-\vartheta|\Delta_n^{-1}\Delta_n^2$ and in turn $\E\big[\big|\sum_{i\in\cI_n} I_i\big|\big] \lesssim |\vartheta_n-\vartheta|$. For the term $II_i$ we use the convergence of $\widehat{\sigmaL{}}_{n,t_{i-1}}$ and the conditioning argument in a similar way as in Step 1 to obtain $\sum_{i\in\cI_n} II_i \rightarrow 0$ in probability. For the term $III_i$, an analysis of the convergence of  $\widehat{\sigmaL{}}_{n,t_{i-1}}$ using Assumption \ref{asm:sigmasHolder} shows that $\E[\sum_{i \in \mathcal I_n}|III_i|] \lesssim \Delta_n^{\alpha}$.
This proves 
$$\Delta_n\sum_{i\in\cI_n} \big(\widetilde{\ell}_{\vartheta_n, \widehat{\sigmaL{}}_n}^{i,n})^2-\Delta_n\sum_{i\in\cI_n} \big(\widetilde{\ell}_{\vartheta,\sigmaL{}}^{i,n}\big)^2 \rightarrow 0$$ in probability and the result remains true with $\widehat{\vartheta}_{2,n}$ in place of $\vartheta_n$. Since
$\Delta_n\sum_{i\in\cI_n} \big(\widetilde{\ell}_{\vartheta,\sigmaL{}}^{i,n}\big)^2 \sim \widetilde{I}_{\vartheta,\sigmaS{},\sigmaL{}}$ in probability under $\mathbb P_{\vartheta,\sigmaS{},\sigmaL{}}$ 
we obtain \eqref{eqn:convergenceVersInformationAvecEstimateursSigma}.\\

\noindent {\bf Step 4.} By definition of $\widetilde{\vartheta}_{2,n}$, we have
\begin{align*}
        \Delta_n^{-1/2}\big(\widetilde{\vartheta}_{2,n}-\vartheta\big)\widetilde{I}_{\vartheta,\sigmaS{},\sigmaL{}} =&\Delta_n^{-1/2}(\widehat{\vartheta}_{2,n}-\vartheta)\widetilde{I}_{\vartheta,\sigmaS{},\sigmaL{}}+\Delta_n^{-1/2}\frac{\Delta_n\widetilde{I}_{\vartheta,\sigmaS{},\sigmaL{}}\sum_{i\in\cI_n}\widetilde{\ell}_{\widehat \vartheta_{2,n}, \widehat{\sigmaL{}}_n}^{i,n} }{\Delta_n\sum_{i\in\cI_n} \big(\widetilde{\ell}_{\widehat \vartheta_{2,n}, \widehat{\sigmaL{}}_n}^{i,n}\big)^2}
        \sim \Delta_n^{1/2}\sum_{i\in\cI_n} \widetilde{\ell}_{\vartheta,\sigmaL{}}^{i,n}
\end{align*}
in probability under $\mathbb P_{\vartheta,\sigmaS{},\sigmaL{}}$ thanks to \eqref{eqn:dlAnA} and \eqref{eqn:convergenceVersInformationAvecEstimateursSigma} established in the two previous steps. We further write
$$\Delta_n^{1/2}\sum_{i\in\cI_n} \widetilde{\ell}_{\vartheta,\sigmaL{}}^{i,n} = \Delta_n^{1/2}\sum_{i=1}^n \widetilde{\ell}_{\vartheta,\sigmaL{}}^{i,n} - \Delta_n^{1/2}\sum_{i=1}^{\lfloor h_n\Delta_n^{-1} \rfloor} \widetilde{\ell}_{\vartheta,\sigmaL{}}^{i,n},$$
and we claim that 
\begin{equation} \label{eq: le dernier reste}
\Delta_n^{1/2}\sum_{i=1}^{\lfloor h_n\Delta_n^{-1} \rfloor} \widetilde{\ell}_{\vartheta,\sigmaL{}}^{i,n} \rightarrow 0
\end{equation}
in probability. We conclude the proof using the following limit theorem, proof of which is delayed until Appendix \ref{proof of TCL}.
\begin{lem} \label{lem:convergenceFctScoreEfficace}
       Work under Assumptions~\ref{asm:regulariteFonctions} and \ref{asm:sigmasHolder} with $\alpha>1/2$. Then
        $$\Delta_n^{1/2}\sum_{i=1}^n\widetilde{\ell}_{\vartheta,\sigmaL{}}^{i,n} \rightarrow \mathcal N(0, \widetilde{I}_{\vartheta,\sigmaS{},\sigmaL{}})$$
        stably in law, where, conditional on $\mathcal F$, the random variable $\mathcal N(0, \widetilde{I}_{\vartheta,\sigmaS{},\sigmaL{}})$ is centred Gaussian with conditional variance $\widetilde{I}_{\vartheta,\sigmaS{},\sigmaL{}}$.
\end{lem}
In view of Lemma~\ref{lem:convergenceFctScoreEfficace} we obtain Theorem \ref{thm:estimationEfficace}.\\

\noindent {\bf Step 5.} It remains to prove \eqref{eq: le dernier reste}.
Write $\Delta_n^{1/2}\widetilde{\ell}_{\vartheta, \sigmaS{},\sigmaL{}}^i=\frac{(T_2-T_1)e^{-\vartheta (T_2-T_1)}}{(1-e^{-\vartheta(T_2-T_1)})^3}(I_i+II_i)$, with
\begin{align*}
        I_i &=\Delta_n^{1/2}\tfrac{(\Delta_i^n X^2-\Delta_i^n X^1)(\Delta_i^n X^2-e^{-\vartheta(T_2-T_1)}\Delta_i^n X^1)}{\Delta_n \sigmaL{t_{i-1}}^2},\\
        II_i &=\Delta_n^{1/2}(\Delta_i^n X^2-\Delta_i^n X^1)(\Delta_i^n X^2-e^{-\vartheta(T_2-T_1)}\Delta_i^n X^1)\tfrac{\int_{t_{i-1}}^{t_i}(\sigmaL{t_{i-1}}^2-\sigmaL{t}^2)dt}{\Delta_n\sigmaL{t_{i-1}}^2\int_{t_{i-1}}^{t_i}\sigmaL{t}^2dt}.
\end{align*}
We readily have
$$\E\big[I_i^2]\leq \tilde{c}^{-2}\Delta_n^{-1}\E\big((\Delta_i^n X^2-\Delta_i^n X^1)^2(\Delta_i^n X^2-e^{-\vartheta(T_2-T_1)}\Delta_i^n X^1)^2\big)\lesssim \Delta_n,$$
so that $\E[(\sum_{i=1}^{\lfloor h_n\Delta_n^{-1} \rfloor} I_i)^2] = \sum_{i=1}^{\lfloor h_n\Delta_n^{-1} \rfloor}\E[I_i^2] \lesssim h_n\Delta_n^{-1}\Delta_n=h_n \rightarrow 0$. The term $II_i$ is similar to the term $I$ in Step 1. We readily obtain 
obtain 
$\E[|II_i|]\lesssim \Delta_n^{1/2+\alpha}$. It follows that 
\[
        \E\big[\big|\sum_{i=1}^{\lfloor h_n\Delta_n^{-1} \rfloor} II_i \big|\big] \lesssim h_n \Delta_n^{-1}\Delta_n^{1/2+\alpha}=h_n \Delta_n^{\alpha-1/2} \rightarrow 0
\]
since $\alpha\geq 1/2$, and \eqref{eq: le dernier reste} follows. The proof of Theorem \ref{thm:estimationEfficace} is complete.

\section{Appendix}

\subsection{Proof of Lemma~\ref{technical lemma 1}} \label{subsection:preuveLemmeTechnique1}

\subsubsection*{The first part of the result} Since  $\sup_{t \in [0,T]} t^{-s}\omega(Y)_t <\infty$ for some $s>1/2$, we have that $Y$ is continuous in probability on $[0,T]$.
Write
\[
        \Delta_n^{-1}\sum_{i=1}^n\big(\overline{\Delta}_i^n(Y)\big)^2-\int_0^TY_t^2dt=I+II+III,
\]
with
\begin{align*}
        I &=\sum_{i=1}^n\int_{(i-1)\Delta_n}^{i\Delta_n} (Y_{(i-1)\Delta_n}^2-Y_t^2)dt,\\
        II &=\Delta_n^{-1}\sum_{i=1}^n \Big(\int_{(i-1)\Delta_n}^{i\Delta_n}(Y_t-Y_{(i-1)\Delta_n})dt\Big)^2,\\
        III &=2\sum_{i=1}^n Y_{t_{i-1}}\int_{(i-1)\Delta_n}^{i\Delta_n}(Y_t-Y_{(i-1)\Delta_n})dt.
\end{align*}
First, fix $\epsilon>0$. There exists some $\eta>0$ such that $\E[|Y_t-Y_s|]<\epsilon$ as soon as $|t-s|<\eta$. Moreover, by localisation we may (and will) assume that there is some $M>0$ such that $\sup_t|Y_t|\leq M$. It follows that
\[
        \E[|I|]\leq \sum_{i=1}^n\int_{(i-1)\Delta_n}^{i\Delta_n} \E[|Y_{(i-1)\Delta_n}^2-Y_t^2|]dt \leq 2 M  \sum_{i=1}^n\int_{(i-1)\Delta_n}^{i\Delta_n} \E[|Y_{(i-1)\Delta_n}-Y_t|]dt \leq 2TM\varepsilon
\]
as soon as $\Delta \leq \eta$ which is true for large enough $n$.
Thus $I \rightarrow 0$ in probability. The proof is similar for $II$ and $III$.
\subsubsection*{The second part of the result} Write
$$
\Delta_n^{-1}\sum_{i=1}^n \overline{\Delta}_i^n(Y) \int_{(i-1)\Delta_n}^{i\Delta_n} Z_td\Ws{t} - \sum_{i=1}^n \int_{(i-1)\Delta_n}^{i\Delta_n} Y_t Z_t d\Ws{t} 
      = I+II,
      $$
      with
      \begin{align*}
        I & = \sum_{i=1}^n \int_{(i-1)\Delta_n}^{i\Delta_n} (Y_{(i-1)\Delta_n}-Y_t) Z_t d\Ws{t} \\
        II & = \Delta_n^{-1}\sum_{i=1}^n \big(\overline{\Delta}_i^n(Y)-\Delta_n Y_{(i-1)\Delta_n}\big) \int_{(i-1)\Delta_n}^{i\Delta_n} Z_td\Ws{t}. 
\end{align*}
Fix $\epsilon >0$ and $\eta>0$ such that $\E[|Y_t-Y_s|]<\epsilon$ as soon as $|t-s|<\eta$. By localisation, we may assume that $Z$ is such that $\sup_t\max(|Z_t|, |Y_t|) \leq M$. By the martingale property, 
$$
      \E[I^2]
        =\sum_{i=1}^n\E\big[\big(\int_{(i-1)\Delta_n}^{i\Delta_n}(Y_{(i-1)\Delta_n}-Y_t) Z_t d\Ws{t}\big)^2\big] 
        = \sum_{i=1}^n\int_{(i-1)\Delta_n}^{i\Delta_n}\E\big[(Y_{(i-1)\Delta_n}-Y_t)^2 Z_t^2\big]dt \leq 2TM^3\epsilon
        $$
 as soon as $\Delta_n \leq \eta$ which is true for large enough $n$.
For $II$, by Cauchy-Schwarz inequality,
\[
       II \leq \big(\sum_{i=1}^n \big(\Delta_n^{-1}\int_{t_{i-1}}^{t_i}(Y_t -Y_{t_{i-1}})dt\big)^2\big)^{1/2} \big(\sum_{i=1}^n \big(\int_{t_{i-1}}^{t_i} Z_td\Ws{t}\big)^2\Big)^{1/2} \lesssim \big(\sum_{i=1}^n \big(\Delta_n^{-1}\int_{t_{i-1}}^{t_i}(Y_t -Y_{t_{i-1}})dt\big)^2\big)^{1/2}
\]
in probability.
Let $\varphi(x)=\indiq_{[0,1)}$ be the Haar function, and $\varphi_{j,k}(x)=2^{j/2}\varphi(2^jx-k)=2^{j/2}\indiq_{[k2^{-j},(k+1)2^{-j})}$ for any $j\geq 0,k\in\Z$. We prove the result under the restriction that $n=2^j$ and that the $t_i$ are of the form $k2^{-j}$. The general case of a regular mesh $t_i = i\Delta_n$ is slightly more intricate but follows the same ideas. We have
$$
        \sum_{i=1}^n \big(\Delta_n^{-1}\int_{t_{i-1}}^{t_i}(Y_t -Y_{t_{i-1}})dt\big)^2 =\frac{1}{T^2}\sum_{k=1}^{2^j} \big(2^j\int_{t_{k-1}}^{t_k}(Y_t -Y_{t_{k-1}})dt\big)^2
       =\sum_{k=0}^{2^j-1} \big((P_j(Y_{\cdot T})(k2^{-j}) - Y_{k2^{-j}T}\big)^2
$$
where $P_j(f)=\sum_{k=0}^{2^j-1}(\int\varphi_{j,k}f)\varphi_{j,k}$ is the orthogonal projection on $\mathrm{Span}\{\varphi_{j,k},k=1,\ldots, 2^j\}$. For large enough $j$, there exists some constant $C>0$ such that
\[
        \E\big[\sum_{k=0}^{2^j-1} \big(P_j(Y_{\cdot T})\big(\frac{k}{2^j}\big) - Y_{\frac{k}{2^j}T}\big)^2\big]\leq C 2^j\E\big[\int_0^1 (P_j(Y_{\cdot T})(u)-Y_{u T})^2du\big].
\]
It follows that
\begin{align*}
        \E\big[\sum_{i=1}^n \big(\Delta_n^{-1}\int_{t_{i-1}}^{t_i}(Y_t -Y_{t_{i-1}})dt\big)^2\big]&\leq C 2^j \E\big[\int_0^1 (P_j(Y_{\cdot T})(u)-Y_{u T})^2du\big] \\
        &\lesssim 2^j \E\big[\int_0^1 \big(\int_0^1 2^j \indiq_{|u-y|\in[0,2^{-j})}(Y_{yT}-Y_{uT})dy\big)^2du\big] \\
         & \lesssim   2^j\big(\int_0^1\omega_{2^{-j}xT}(Y)dx\big)^2 
         \lesssim 2^{j(1-2s)} 
\end{align*}
and this term converges to $0$ since $s>1/2$ by assumption. The proof of Lemma~\ref{technical lemma 1} is complete.

\subsection{Proof of Lemma \ref{lem:convergenceFctScoreEfficace}} \label{proof of TCL}
First we establish the result for
\[
        \frac{T_2-T_1}{e^{\vartheta(T_2-T_1)}-1}\Delta_n^{1/2}\sum_{i=1}^n\chi_i^n,\;\;\text{with}\;\;
        \chi_i^n=\frac{\sigmaS{t_{i-1}}\Delta_i^n\Wl{} \int_{t_{i-1}}^{t_i}e^{-\vartheta(T_1-t)}d\Ws{t}}{\Delta_n \sigmaL{(i-1)\Delta_n}},
\]
by applying Lemma~3.7 in \cite{bib:statisticsAndHFDataInStatisticalMethodsForSDE}. To do so, we check conditions~(3.43)--(3.46) in \cite{bib:statisticsAndHFDataInStatisticalMethodsForSDE} for $\chi_i^n$. We keep up with the notation of  \cite{bib:statisticsAndHFDataInStatisticalMethodsForSDE}. 
First, we have $\E[\chi_i^n|\cF_{i-1}]=0$, which ensures (3.43) with $A_t=0$.
Next,
\begin{align*}
        \E[(\chi_i^n)^2|\cF_{i-1}]&=\Delta_n\tfrac{\sigmaS{(i-1)\Delta_n}^2}{\Delta_n^2\sigmaL{(i-1)\Delta_n}^2}\E[(\Delta_i^n \Wl{})^2|\cF_{i-1}]\E[(\int_{(i-1)\Delta_n}^{i\Delta_n}e^{-\vartheta(T_1-t)}d\Ws{t})^2|\cF_{i-1}] \\
 &       =\tfrac{\sigmaS{(i-1)\Delta_n}^2}{\sigmaL{(i-1)\Delta_n}^2}\int_{(i-1)\Delta_n}^{i\Delta_n}e^{-2\vartheta(T_1-t)}dt,
\end{align*}
so that
\[
        \sum_{i=1}^{\lfloor t/\Delta_n \rfloor} \E[(\chi_i^n)^2|\cF_{i-1}] \to \int_0^t \tfrac{e^{-2\vartheta(T_1-s)}\sigmaS{s}^2}{\sigmaL{s}^2}ds
\]
in probability. This is Condition (3.44) in \cite{bib:statisticsAndHFDataInStatisticalMethodsForSDE} with $C_t=\int_0^t \frac{e^{-2\vartheta(T_1-s)}\sigmaS{s}^2}{\sigmaL{s}^2}ds$. It follows that 
\begin{align*}
        \E[(\chi_i^n)^4|\cF_{i-1}]&=\Delta_n^2\tfrac{\sigmaS{(i-1)\Delta_n}^4}{\Delta_n^4\sigmaL{(i-1)\Delta_n}^4}\E[(\Delta_i^n \Wl{})^4|\cF_{i-1})\E[(\int_{(i-1)\Delta_n}^{i\Delta_n}e^{-\vartheta(T_1-t)}d\Ws{t})^4|\cF_{i-1}]\\
        &\leq 9\Delta_n^2\tfrac{\sigmaS{(i-1)\Delta_n}^4}{\sigmaL{(i-1)\Delta_n}^4}\Big(\int_{(i-1)\Delta_n}^{i\Delta_n}e^{-2\vartheta(T_1-t)}dt\Big)^2
\end{align*}
by independence of the two Wiener integrals. Therefore $\sum_{i=1}^n\E[(\chi_i^n)^4|\cF_{i-1}] \rightarrow 0$ in probability. This is condition~(3.45) in \cite{bib:statisticsAndHFDataInStatisticalMethodsForSDE}. Finally, 
 $\E\big[\chi_i^n\Delta_i^n B\big|\cF_{i-1}\big]= \E\big[\chi_i^n\Delta_i^n \overline{B}\big|\cF_{i-1}\big] = 0$ by independence which ensures condition~(3.46)  in \cite{bib:statisticsAndHFDataInStatisticalMethodsForSDE}. We subsequently apply 
Lemma~3.7 in \cite{bib:statisticsAndHFDataInStatisticalMethodsForSDE} to conclude that
\[
        \frac{T_2-T_1}{e^{\vartheta(T_2-T_1)}-1}\Delta_n^{1/2}\sum_{i=1}^n\chi_i^n
\]
converges stably in law to a random variable, which, conditional on $\cF$, is Gaussian with variance $\widetilde{I}_{\vartheta,\sigmaS{},\sigmaL{}}$.
In order to complete the proof, we write
\[
        \Delta_n^{1/2}\sum_{i=1}^n \widetilde{\ell}_{\vartheta,\sigmaL{}}^{i,n} = \Delta_n^{1/2} \sum_{i=1}^n \Big(\widetilde{\ell}_{\vartheta,\sigmaL{}}^{i,n}-\frac{T_2-T_1}{e^{\vartheta(T_2-T_1)}-1}\sum_{i=1}^n\chi_i^n\Big)+\frac{T_2-T_1}{e^{\vartheta(T_2-T_1)}-1}\Delta_n^{1/2}\sum_{i=1}^n\chi_i^n
\]
and it remains to show the convergence of $\Delta_n^{1/2}\sum_{i=1}^n \big(\widetilde{\ell}_{\vartheta,\sigmaL{}}^{i,n}-\frac{T_2-T_1}{e^{\vartheta(T_2-T_1)}-1}\sum_{i=1}^n\chi_i^n\big)\rightarrow 0$ in probability. This is done using similar arguments as used in the proof of Theorem~\ref{thm:estimationNPunePeriode}. We omit the details.\\

\noindent {\bf Acknowledgements.} {\it We are grateful to the suggestions and comments of two anonymous referees that helped to considerably improve a former version of the manuscript.}

\bibliographystyle{plain}
\bibliography{Biblio}

\end{document}